\documentclass[10pt,a4paper]{article}
\usepackage{mathrsfs}
\usepackage{amsfonts}
\begin{document}

\begin{center}
 {\huge\textbf{Solutions of Cauchy problem for multiple inhomogeneous wave equation}}
\vspace{0.3cm}

\begin{center}
\rm Guang-Qing Bi
\end{center}

\end{center}

\noindent
\begin{abstract} We define a class of pseudo-differential operators in a completely new way, which is called the abstract operators and expounded systematically the theory of abstract operators. By combining abstract operators with the Laplace transform, we can apply the Laplace transform to any $n+1$ dimensional linear higher-order partial differential equations $P(\partial_x,\partial_t)u=f(x,t)$ directly, without using the Fourier transform. By making introduction of abstract operators $G(\partial_x,t):=\mathcal{L}^{-1}[1/P(\partial_x,s)]$, the analytic solutions of initial value problems are expressed in these abstract operators, including the multiple inhomogeneous wave equation associated with the shifted Laplace-Beltrami operator on real hyperbolic spaces. By writing abstract operators in this class into integral forms, the solutions in operator form are represented into integral forms. Thus the analytic solutions of Cauchy problem for the multiple wave equation on $\mathbb{R}^n$ can be represented in the integrations of some given functions, without using the traditional Fourier transform technique. As a further application, we study the solvability of initial-boundary value problem for the linear higher-order partial differential equations and deduce new distinguishable method associated with the second-order linear self-adjoint elliptic operators.
\end{abstract}

\footnotetext{\hspace*{-.45cm}\footnotesize
\noindent{\bf 2010 \emph{Mathematics Subject Classification}.}\ \ Primary 35L05, 35S10, 35S11; Secondary 35C15, 44A10, 42A24.\\
\noindent{\bf \emph{Key words and phrases}.}\ \ Abstract operators, Pseudo-differential operators, Wave equation, linear higher-order partial differential equations, Initial value problems, Initial-boundary value problems, Laplace transform.}

\section{Introduction}
\noindent
\renewcommand\theequation{1.\arabic{equation}}

In 1960s, the general theory of linear partial differential equations made important progress by using the generalized function and its Fourier transform, further, the theory of pseudo-differential operators. Building on the theory of Fourier transform, the concept of pseudo-differential operators is defined as follows (See, e.g., Chen \cite[pp. 36-38]{chen}):

Let $a(x,\xi)\in{C}^\infty(\mathbb{R}_x^n\times\mathbb{R}_\xi^n)$. If for any $\alpha,\beta\in\mathbb{N}^n$ and an real number $m$,
\[|\partial_\xi^\alpha\partial_x^\beta{a}(x,\xi)|\leq{C}_{\alpha,\beta}(1+|\xi|)^{m-|\alpha|}\]
is tenable (denoted by $a\in S^m$), where $C_{\alpha,\beta}$ is a constant, then the linear continuous mapping $A$ :
$\mathscr{S}(\mathbb{R}^n)\rightarrow\mathscr{S}(\mathbb{R}^n)$ can be defined as
\begin{equation}\label{chen1}
Au(x):=(2\pi)^{-n}\int_{\mathbb{R}^n}\!\!e^{i\xi{x}}a(x,\xi)\hat{u}(\xi)d\xi,
\end{equation}
which is called the pseudo-differential operators, denoted by
$a(x,D)$, where $a(x,\xi)$ is the symbols of pseudo-differential operators $a(x,D)$.

Gala consider the regularity of the following type (See \cite[p. 1262]{gala}):
\[|\partial_\xi^\alpha{a}(x,\xi)|\leq{A}_\alpha(1+|\xi|)^{-|\alpha|}\quad(\alpha\in\mathbb{N}^n),\]
where the variable $x$ is related to a suitable positive function $\omega$ such that
\[|\partial_\xi^\alpha{a}(x+y,\xi)-\partial_\xi^\alpha{a}(x,\xi)|\leq{A}_\alpha\omega(|y|)(1+|\xi|)^{-|\alpha|}.\]
Lannes \cite[p. 495]{lann} consider pseudo-differential operators whose symbol $a(x,\xi)$ is not even infinitely smooth with respect to $x$.
Therefore, we can introduce more general classes of pseudo-differential operators by relaxing the restrictions on symbols according to the problem being discussed.

In fact, the biggest constraint on the symbols should be derived from the application of Fourier transform in the definition of pseudo-differential operators.
However, the definition of pseudo-differential operators means that the following relationship is tenable:
\begin{equation}\label{chen2}
a(x,\xi):=e^{-i\xi{x}}a(x,D)e^{i\xi{x}}\quad\mbox{or}\quad a(x,D)e^{i\xi{x}}=a(x,\xi)e^{i\xi{x}},
\end{equation}
where $\xi{x}:=\langle\xi,x\rangle=\xi_1x_1+\xi_2x_2+\cdots+\xi_nx_n$. Obviously, the properties of $a(x,D)$ depend on $a(x,\xi)$. equality
(\ref{chen2}) indicates the mapping relation between $a(x,D)$ and $a(x,\xi)$. For a given function $a(x,\xi)$, if the corresponding algorithms of $a(x,D)$ can be derived from this type of mapping relations, to determine the domain and range of $a(x,D)$, then this type of mapping relations can be the best definition of pseudo-differential operators. Luckily, this idea totally works by making use of the Taylor theorem of analytic function and analytic continuation technique. Based on this idea, the author introduce the mathematical concept of abstract operators, including  $\exp(h\partial_x),\;\sin(h\partial_x),\;\cos(h\partial_x),\;\sinh(h\partial_x)\;\mbox{and}\;\cosh(h\partial_x)$,
where $h\partial_x:=\langle{h,\partial_x}\rangle=h_1\partial_{x_1}+h_2\partial_{x_2}+\cdots+h_n\partial_{x_n}$, as well as
\[\exp(tP(\partial_x)),\;\cos(tP(\partial_x)^{1/2}),\;\frac{\sin(tP(\partial_x)^{1/2})}{P(\partial_x)^{1/2}},\;
\cosh(tP(\partial_x)^{1/2})\;\mbox{and}\;\frac{\sinh(tP(\partial_x)^{1/2})}{P(\partial_x)^{1/2}}\] derived by the abstract operators $G(P(\partial_x),t):=\mathcal{L}^{-1}[1/P(\partial_x,s)]$ for any $n+1$ dimensional linear partial differential equations $P(\partial_x,\partial_t)u=f(x,t)$. This means that one can solve linear partial differential equations explicitly by using the theory of abstract operators. Thus the general solving procedure of initial value problem for any linear higher-order partial differential equations is derived clearly in a different way than the Fourier transform by Guang-Qing Bi \cite{bi97} to \cite{bi11}. In this paper, we will obtain some new important solving formulas, and further discuss the solvability of initial-boundary value problem for linear higher-order partial differential equations. However, the concept of abstract operators is not fully expressed since it in a minimalist form first introduced in 1997. Therefore, our primary task in this paper is to expound systematically the fundamental theory of abstract operators as well as their applications in linear higher-order partial differential equations including the multiple inhomogeneous wave equation on $\mathbb{R}^n$ or $\mathbb{H}^n$. The theory and methodology of abstract operators are not only used to establish the theory of linear partial differential equations, but also are gradually integrating with other mathematical branches to obtain more extensive applications of pseudo-differential operators.

\section{Preliminaries (Theory of abstract operators)}
\noindent\setcounter{equation}{0}
\renewcommand\theequation{2.\arabic{equation}}

In this section we will see how the Analytic Continuity Fundamental Theorem can be applied in defining pseudo-differential operators without using Fourier transform.

\subsection{Basic notions}
\noindent

For $\alpha:=(\alpha_1,\cdots,\alpha_n)\in\mathbb{R}^n$, if $\alpha_i\geq0$, then $\alpha$ is called a multiple indicator, denoted by $\alpha\in\mathbb{N}^n$. For $\alpha:=(\alpha_1,\cdots,\alpha_n)\in\mathbb{N}^n$, we have universal mark as follows
\[\partial_x^\alpha u(x):=\frac{\partial^{|\alpha|}u(x)}{\partial x_1^{\alpha_1}\partial x_2^{\alpha_2}\cdots\partial x_n^{\alpha_n}}:=\partial_{x_1}^{\alpha_1}\partial_{x_2}^{\alpha_2}\cdots\partial_{x_n}^{\alpha_n}u,
\quad|\alpha|:=\sum^n_{i=1}\alpha_i,\quad\partial_{x_i}:=\partial_i.\]

\textbf{Definition 2.1.} Within each convergence circle of analytic functions, if the effects on the series term by term from the linear operator converge uniformly to the effects on the sum function, then the operator is called having the analytic continuity.

\textbf{Definition 2.2.} Let $\alpha\in\mathbb{N}^n,\,x\in{\mathbb{R}^n}$. $x^\alpha$ can be called the base functions of analytic function, which with exponential form $e^{\xi x}$ can be defined by the following  eigenfunctions of $\partial_x$:
\[\left\{\begin{array}{l@{\qquad}l}
\partial_xu(x)=\xi u(x) & \xi\in\mathbb{R}^n,\\
u(0)=1.\end{array}\right.\]
Here $\alpha$ and $\xi$ are called the characters of $x^\alpha$ and $e^{\xi{x}}$ respectively.

In view of the Taylor theorem of analytic function (include binomial formula of integer power)
\[f(z)=\sum^\infty_{|\alpha|=0}\left.\partial_z^\alpha f(z)\right|_{z=s}\frac{(z-s)^\alpha}{\alpha!},\quad\alpha\in\mathbb{N}^n,\;z\in\mathbb{C}^n\]
and using analytic continuation technique, we can use the Definition 2.1 and Definition 2.2 to deduce the following fundamental principle:

\textbf{Analytic Continuity Fundamental Theorem.} Let $A$ and $B$ be the linear operators with the analytic continuity, their domains be denoted by $\mathscr{D}(A)$ and $\mathscr{D}(B)$ respectively. If the variable $x\in\mathbb{R}^n$ is related to $y\in\mathbb{R}^n$ such that one of the following two operator identities
\[Ax^\alpha=By^\alpha\quad\mbox{or}\quad Ae^{\xi{x}}=Be^{\xi y}\]
tenable, where the expressions of $A$ and $B$ do not explicitly contain the characters $\alpha\in\mathbb{N}^n$ or $\xi\in\mathbb{R}^n$, then
\[Af(x)=Bf(y),\qquad\forall f(x),f(y)\in\mathscr{D}(A)\cap\mathscr{D}(B),\]
which is a slightly modified version of the result given by Guang-Qing Bi \cite[pp. 7-8]{bi97}.

Building on the Analytic Continuity Fundamental Theorem, we now develop the following mathematical concept of abstract operators to define a class of pseudo-differential operators in a new way:

\textbf{Definition 2.3.} Let $T>0,\,\Omega\subseteq\mathbb{R}^n$ be an open set, a class of linear operators with the analytic continuity is
called the abstract operators taking
$x\partial_x:=(x_1\partial_{x_1},\cdots,x_n\partial_{x_n})=(x_1\partial_1,\cdots,x_n\partial_n)$ as the operator element, denoted by
$f(t,x\partial_x)$,
$x\in\Omega,\,t\in(0,T]$, which is also the linear continuous mapping $f(t,x\partial_x)$:
$C^\infty(\Omega)\rightarrow{C^\infty(\Omega)}$, and if it acts on the base functions $x^\alpha$, we have
\begin{equation}\label{d2.3}
  f(t,\,x\partial_x)x^\alpha:=f(t,\alpha)x^\alpha,\quad\forall{f(t,\alpha)}\in{C^\infty}((0,T]\times\mathbb{R}^n),
\end{equation}
where $f(t,\alpha),\alpha\in\mathbb{N}^n$ is called the symbols of abstract operators $f(t,x\partial_x)$.

\textbf{Definition 2.4.} A class of linear operators with the analytic continuity is called the abstract operators taking $\partial_x:=(\partial_{x_1},\cdots,\partial_{x_n})$ as the operator element, denoted by $f(t,\partial_x),\,x\in\Omega,\,t\in(0,T]$, which is also the linear continuous mapping $f(t,\partial_x)$:
$C^\infty(\Omega)\rightarrow{C^\infty(\Omega)}$, and if it acts on the exponential base functions $e^{\xi{x}}$, we have
\begin{equation}\label{d2.4}
  f(t,\partial_x)e^{\xi{x}}:=f(t,\xi)e^{\xi{x}},\quad\forall{f(t,\xi)}\in{C^\infty}((0,T]\times\mathbb{R}^n),
\end{equation}
where $f(t,\xi),\,\xi\in\mathbb{R}^n$ is called the symbols of abstract operators $f(t,\partial_x)$.

\textbf{Remark 2.1.} By making use of (\ref{chen1}) and (\ref{d2.4}), we have
\begin{equation}\label{cd}
  f(t,\partial_x)u(x)=(2\pi)^{-n}\int_{\mathbb{R}^n}\!\!e^{i\xi{x}}f(t,i\xi)\hat{u}(\xi)d\xi,\quad\forall u(x)\in\mathscr{S}(\mathbb{R}^n),
\end{equation}
when $f(t,i\xi)\in S^m$. Therefore, the pseudo-differential operators in view of the Fourier transform can also be called the abstract operators defined on $\mathscr{S}(\mathbb{R}^n)$. Conversely, the abstract operators in view of the Analytic Continuity Fundamental Theorem can also be called the pseudo-differential operators defined on $C^\infty(\Omega)$.

\textbf{Definition 2.5.} The operator identities that determine the domain and range of abstract operators are called the algorithms of the abstract operators.

\textbf{Definition 2.6.} The relational expression between each component of the characters
$\alpha=(\alpha_1,\alpha_2,\ldots,\alpha_n)$ or $\xi=(\xi_1,\xi_2,\ldots,\xi_n)$ is called the characteristic equation.

\textbf{Definition 2.7.} Let $A$ be a linear operator having the analytic continuity, if there exists another linear operator, denoted by $A^{-1}$ such that $AA^{-1}=A^{-1}A=I$, then $A^{-1}$ is called the inverse operator of $A$.

By Definition 2.4 and the Analytic Continuity Fundamental Theorem we obtain:

\textbf{Corollary 2.1.} The operator algebras formed by all abstract operators $f(t,\partial_x)$, are isomorphic to the algebras formed by their symbols $f(t,\xi)$. This isomorphism is determined by $f(t,\partial_x)\leftrightarrow{f}(t,\xi)$, and

\parbox{12cm}{\begin{eqnarray*}\label{c2.1}
f(t,\partial_x)\pm{g}(t,\partial_x) &\leftrightarrow& f(t,\xi)\pm{g}(t,\xi),\\
f(t,\partial_x)\circ{g}(t,\partial_x) &\leftrightarrow& {f}(t,\xi)g(t,\xi).
              \end{eqnarray*}}\hfill\parbox{2cm}{\begin{eqnarray}\end{eqnarray}}

\textbf{Remark 2.2.} It is easily seen from the Corollary 2.1 that \[\cos(ix\partial_y)=\cosh(x\partial_y)\quad\mbox{and}\quad\sin(ix\partial_y)=i\sinh(x\partial_y),\]
where $x\partial_y:=\langle{x,\partial_y}\rangle=x_1\partial_{y_1}+\cdots+x_n\partial_{y_n}$.
Especially, the abstract operators $f(\partial_x)$ and $g(\partial_x),\;x\in\mathbb{R}^n$ are each other's inverse operators, if and only if their symbols $f(\xi)$ and $g(\xi)$ satisfy the algebraic relationship $f(\xi)g(\xi)=1,\;\xi\in\mathbb{R}^n$.

By Definition 2.4 we have

\textbf{Corollary 2.2.} Let $g(x)\in{C^\infty(\mathbb{R}^n)}$. If $g(\partial_\xi)(e^{\xi{x}}f(\xi)),\,\xi\in\mathbb{R}^n$ is continuous at $\xi=\xi_0$, then
\begin{equation}\label{e0}
f(\partial_x)(e^{\xi_0x}g(x))=g(\partial_\xi)(e^{\xi{x}}f(\xi))|_{\xi=\xi_0}.
\end{equation}

By the Analytic Continuity Fundamental Theorem we obtain easily the following results:

\textbf{Example 2.1.} Let $\xi+\partial_x:=(\xi_1+\partial_{x_1},\cdots,\xi_n+\partial_{x_n})$,  $\partial_x+\partial_y:=(\partial_{x_1}+\partial_{y_1},\cdots,\partial_{x_n}+\partial_{y_n})$, $y+x:=(y_1+x_1,\cdots,y_n+x_n)$. Then we have
\begin{equation}\label{e1}
f(\partial_x)(e^{\xi{x}}g(x))=e^{\xi{x}}f(\xi+\partial_x)g(x),\quad\forall g(x)\in C^\infty(\Omega),\;\Omega\in\mathbb{R}^n
\end{equation}
and
\begin{equation}\label{e3}
f(\partial_x+\partial_y)g(x,y)=e^{-x\partial_y}(f(\partial_x)e^{x\partial_y}g(x,y)),\quad\forall g(x,y)\in C^\infty(\Omega),\;\Omega\in\mathbb{R}^n\times\mathbb{R}^n,
\end{equation}
respectively, where $e^{x\partial_y}g(x,y)=g(x,y+x)$.

Taking (\ref{e1}) as the characteristic equation, we can use it and the base functions $e^{\xi{y}}$ to make the following operator equality:
\[f(\partial_x)((e^{\xi{y}}g(x))|_{y=x})=f(\partial_y+\partial_x)(e^{\xi{y}}g(x))|_{y=x}\qquad(y\in\mathbb{R}^n,\;\xi\in\mathbb{R}^n).\]
According to the Analytic Continuity Fundamental Theorem, we can obtain that

\textbf{Example 2.2.} Let $g(x,y)\in C^\infty(\Omega),\;\Omega\in\mathbb{R}^n\times\mathbb{R}^n$. If $g(x,y)$ is continuous at $y=x$, then $f(\partial_x+\partial_y)g(x,y)$ is continuous at $y=x$ such that
\begin{equation}\label{e3ct}
f(\partial_x+\partial_y)g(x,y)|_{y=x}=f(\partial_x)(g(x,y)|_{y=x}).
\end{equation}
Similarly,
\begin{equation}\label{e3ctb}
f(\partial_x+\partial_y+\partial_z)g(x,y,z)|_{z=y=x}=f(\partial_x)(g(x,y,z)|_{z=y=x}),
\end{equation}
where $\partial_x+\partial_y+\partial_z:=(\partial_{x_1}+\partial_{y_1}+\partial_{z_1},\cdots,\partial_{x_n}+\partial_{y_n}+\partial_{z_n})$.

\textbf{Example 2.3.} Let $f(x)\in L^1(\mathbb{R}^n),\,\lambda>0$. Then we have
\begin{equation}\label{e2}
 e^{\lambda|\partial_x|^2}f(x)=\frac{1}{2^n(\lambda\pi)^{n/2}}\int_{\mathbb{R}^n}f(\eta)\exp\left(-\frac{|\eta-x|^2}{4\lambda}\right)d\eta.
\end{equation}

\textbf{Proof.} If taking the following integral formula
\[e^{\lambda\xi_i^2}=\frac{1}{2\sqrt{\pi}}\int^\infty_{-\infty}e^{-\zeta^2/4+\sqrt{\lambda}\,\xi_i\zeta}d\zeta,\;\forall\lambda>0\]
as its characteristic equation, according to the Analytic Continuity Fundamental Theorem we have
\begin{eqnarray*}
  \exp\left(\lambda\frac{\partial^2}{\partial{x_i^2}}\right)f(x)&=&
  \frac{1}{2\sqrt{\pi}}\int^\infty_{-\infty}e^{-\zeta^2/4}f(x_1,\cdots,x_i+\sqrt{\lambda}\zeta,\cdots,x_n)d\zeta\\
   &=&
   \frac{1}{2\sqrt{\pi\lambda}}\int^\infty_{-\infty}f(x_1,\cdots,\eta_i,\cdots,x_n)\exp\left(-\frac{(\eta_i-x_i)^2}{4\lambda}\right)d\eta_i.
\end{eqnarray*}
Similarly, we have (\ref{e2}). Example 2.3 is proved.

\subsection{Algorithms of partial differential operators}
\noindent

In this section we will prove that the concept of abstract operators is the generalization of partial differential operators. In other words,
by applying Definition 2.3 in the case $f(t,\alpha)=\alpha$ and Definition 2.4 in the case $f(t,\xi)=\xi$, we can define the partial differential operators $x\partial_x$ and $\partial_x$ respectively.

\textbf{Corollary 2.3.} Let $\alpha\in\mathbb{N}^n,\;x\in\mathbb{R}^n,\,\xi\in\mathbb{R}^n$. The partial differential operators $x\partial_x$ and $\partial_x$ are the abstract operators defined by
\[x\partial_x x^\alpha:=\alpha x^\alpha\quad\mbox{and}\quad\partial_x e^{\xi{x}}:=\xi{e}^{\xi{x}},\]
respectively. Therefore, all rules of differentiation can be determined by the Analytic Continuity Fundamental Theorem.

\textbf{Proof.} The major parts of rules of differentiation are the derivative principle of function product, the derivative principle of compound function and the chain rule of multivariate function.

Let $x\in\Omega\subseteq\mathbb{R}^1$. If taking $a,b\in\mathbb{R}$ as the characters of the base functions $e^{ax}$ and $e^{bx}$ respectively, then we can combine $e^{ax}$ and $e^{bx}$ with the following characteristic equation
\[{b^n}=\sum^n_{j=0}(-1)^j{n\choose{j}}a^j(a+b)^{n-j}\quad(n\in\mathbb{N})\]
to make the following operator equality for the base functions $e^{ax}$ and $e^{bx}$:
\[e^{ax}\frac{d^n}{dx^n}e^{bx}=\sum^n_{j=0}(-1)^j{n\choose{j}}\frac{d^{n-j}}{dx^{n-j}}\left(e^{bx}\frac{d^j}{dx^j}e^{ax}\right).\]
According to the Analytic Continuity Fundamental Theorem, we have
\begin{equation}\label{1*}
v\frac{d^nu}{dx^n}=\sum^n_{j=0}(-1)^j{n\choose{j}}\frac{d^{n-j}}{dx^{n-j}}\left(u\frac{d^jv}{dx^j}\right),\quad\forall{v,u}\in{C^n}(\Omega).
\end{equation}

Similarly, we can easily derive the following Leibniz rule:
\begin{equation}\label{1}
\frac{d^n}{dx^n}(vu)=\sum^n_{j=0}{n\choose{j}}\frac{d^jv}{dx^j}\frac{d^{n-j}u}{dx^{n-j}},\quad\forall{v,u}\in{C^n}(\Omega).
\end{equation}

Letting $n=1$ in (\ref{1}), then for $v=f(x)$ and $u=g(x)$, we have
\begin{equation}\label{3}
    \frac{d}{dx}(f(x)g(x))=f(x)\frac{d}{dx}g(x)+g(x)\frac{d}{dx}f(x),\quad\forall{f,g}\in{C}^1(\Omega).
\end{equation}

Let $\varphi(x):=f_1(x)f_2(x)\cdots{f}_n(x)$. Generally, we have
\[\frac{d}{dx}\varphi(x)=\sum^n_{j=1}\prod^n_{i=1\atop
i\neq{j}}f_i(x)\frac{d}{dx}f_j(x),\quad\forall{f_j(x)}\in{C}^1(\Omega),\;j=1,2,\cdots,n.\]

If $f_1(x)=f_2(x)=\cdots=f_n(x)=y=g(x)\in{C}^1(\Omega)$, then for differentiable function $y=g(x)$, we have
\[\frac{d}{dx}y^n=ny^{n-1}\frac{dy}{dx}\quad\mbox{or}\quad\frac{d}{dx}y^n=\frac{dy}{dx}\frac{d}{dy}y^n\quad(n\in\mathbb{N}^1).\]
According to the Analytic Continuity Fundamental Theorem, we have the derivative principle of compound function
$f(y)\in{C^1}(V),\,y\in{V}\subseteq\mathbb{R}^1$:
\begin{equation}\label{4}
\frac{d}{dx}f(g(x))=\frac{dy}{dx}\frac{d}{dy}f(y),\qquad{y}=g(x)\in{C}^1(\Omega),\;x\in\Omega\subseteq\mathbb{R}^1.
\end{equation}

By using (\ref{3}) and (\ref{4}), for differentiable functions $x_1(t)$ and $x_2(t)$, $t\in{I}\subseteq\mathbb{R}^1$, we have
\[\frac{d}{dt}(x^{\alpha_1}_1x^{\alpha_2}_2)=\frac{dx^{\alpha_1}_1}{dx_1}\frac{dx_1}{dt}x^{\alpha_2}_2
+x^{\alpha_1}_1\frac{dx^{\alpha_2}_2}{dx_2}\frac{dx_2}{dt},\qquad(\alpha_1,\alpha_2)\in\mathbb{N}^2.\]
Taking this one as the characteristic equation, we can use it and the base functions $x^\alpha$ to make the following operator equality:
\[\frac{d}{dt}x^\alpha=\frac{dx_1}{dt}\frac{\partial}{\partial{x}_1}x^\alpha+\frac{dx_2}{dt}\frac{\partial}{\partial{x}_2}x^\alpha
\qquad(x(t)\in\mathbb{R}^2,\;t\in{I}\subseteq\mathbb{R}^1,\;\alpha\in\mathbb{N}^2).\]
According to the Analytic Continuity Fundamental Theorem,  we have $\forall{f(x)}\in{C^1}(\Omega),\;\Omega\in\mathbb{R}^2$,
\[\frac{d}{dt}f(x)=\frac{dx_1}{dt}\frac{\partial}{\partial{x}_1}f(x)+\frac{dx_2}{dt}\frac{\partial}{\partial{x}_2}f(x),
\quad x_i(t)\in{C}^1(I),\;i=1,2.\]

Similarly, $\forall{f(x)}\in{C^1}(\Omega),\;\Omega\in\mathbb{R}^n$, $t\in{I}\subseteq\mathbb{R}^1$, we can generally derive
\begin{equation}\label{5}
\frac{d}{dt}f(x)=\frac{dx_1}{dt}\frac{\partial}{\partial{x}_1}f(x)+\frac{dx_2}{dt}\frac{\partial}{\partial{x}_2}f(x)
+\cdots+\frac{dx_n}{dt}\frac{\partial}{\partial{x}_n}f(x),\;\;x_i(t)\in{C}^1(I),\,i=1,\cdots,n.
\end{equation}

If taking $n\in\mathbb{N}^1$ as the characters of base functions $z^n$, then the binomial formula can be expressed as the following characteristic equation:
\[(z+h)^n=\sum^\infty_{j=0}\frac{h^j}{j!}\frac{d^j}{dz^j}z^n\quad(z\in\mathbb{C}^1).\]
According to the Analytic Continuity Fundamental Theorem, the Taylor formula is tenable for any analytic function, namely
\begin{equation}\label{2}
f(z+h)=\sum^\infty_{j=0}\frac{h^j}{j!}\frac{d^j}{dz^j}f(z),\quad\forall f(z)\in{C^\infty}(\Omega),\;\Omega\subseteq\mathbb{C}^1,\;|h|<R.
\end{equation}

Let $a$ and $b$ be real numbers such that $a\leq{b}$. Without losing the universality, assuming $a\geq0,\;b\geq0$ we have
\[na^{n-1}\leq(a^{n-1}+a^{n-2}b+a^{n-3}b^2+\cdots+ab^{n-2}+b^{n-1})\leq{n}b^{n-1},\quad{n}=1,2,\cdots.\]
namely
\[na^{n-1}\leq\frac{b^n-a^n}{b-a}\leq{n}b^{n-1}\quad\mbox{or}\quad
{a}\leq\sqrt[n-1]{\frac{1}{n}\frac{b^n-a^n}{b-a}}\leq{b}.\]
Therefore, if $a\leq{b}$, then we have $a\leq{c}\leq{b}$, making
\[nc^{n-1}=\frac{b^n-a^n}{b-a}\quad\mbox{or}\quad\left.\frac{d}{dx}x^n\right|_{x=c}=\frac{b^n-a^n}{b-a}\qquad(n\in\mathbb{N}^1).\]
Taking this one as the characteristic equation, according to the Analytic Continuity Fundamental Theorem, we obtain the following Lagrange mean value theorem:

If $a\leq{b}$, then there exists $c\in[a,b]$ such that
\begin{equation}\label{lg}
\left.\frac{d}{dx}f(x)\right|_{x=c}=\frac{f(b)-f(a)}{b-a},\quad\forall{f(x)}\in{C^1}[a,b].
\end{equation}
Thus Corollary 2.3 is proved.

\textbf{Remark 2.3.} Constructing the operator equality tenable for the base functions by the characteristic equation through the definition of abstract operators, and then deducing that it is also tenable for any analytic functions according to the Analytic Continuity Fundamental Theorem, thus we derive new operator identities. Finding or establishing a new operator identities requires the knowledge of the corresponding characteristic equation in advance, without knowing the specific form of the new operator identities. Therefore, it all boils down to seek or construct appropriate characteristic equations. The key to transform the characteristic equation to the corresponding operator identities, is constructing the operator equality true for the base functions by using the specific form of the characteristic equation and the definition of abstract operators, and also, only when the operators constructed are all linear and the expressions of the linear operators do not explicitly contain the characters of the base functions can we derive that the operator equality is not only true for the base functions, but also for any analytic functions, according to the Analytic Continuity Fundamental Theorem.

\textbf{Example 2.4.} Let $y=x^2,\,x\in\mathbb{R}^1$. $\forall{f(y)}\in{C^k}(\Omega),\;y\in\Omega\subseteq\mathbb{R}^1$, we have
\begin{equation}\label{y5}
\frac{d^k}{dx^k}f(x^2)=\sum^{[k/2]}_{j=0}\frac{k!}{j!\,(k-2j)!}(2x)^{k-2j}\frac{d^{k-j}}{dy^{k-j}}f(y).
\end{equation}
If $y=x^3,\,x\in\mathbb{R}^1$, then $\forall{f(y)}\in{C^n}(\Omega)$,
\begin{equation}\label{y5.3}
\frac{d^n}{dx^n}f(x^3)=\sum^{[n/2]}_{k=0}\sum^k_{j=0}\frac{3^{n-k-2j}n!}{j!\,(k-j)!\,(n-2k-j)!}x^{2n-3k-3j}\frac{d^{n-k-j}}{dy^{n-k-j}}f(y).
\end{equation}

\textbf{Proof.} Let $x\in\mathbb{R}^1,\,\xi\in\mathbb{R}^1,\,\zeta=3x^2\lambda$, $\lambda$ be a real parameter. By applying Corollary 2.2 and (\ref{e1}) of Example 2.1, we have
\begin{eqnarray*}
  f(\partial_x)e^{\lambda{x^3}} &=& \left.\exp\left(\lambda\frac{\partial^3}{\partial\xi^3}\right)(e^{\xi{x}}f(\xi))\right|_{\xi=0}
  =\left.\exp\left(\lambda\left(x+\frac{\partial}{\partial\xi}\right)^3\right)f(\xi)\right|_{\xi=0}\\
   &=& \left.\exp\left(\lambda{x^3}+3x^2\lambda\frac{\partial}{\partial\xi}+3x\lambda\frac{\partial^2}{\partial\xi^2}+\lambda\frac{\partial^3}{\partial\xi^3}\right)f(\xi)\right|_{\xi=0}\\
   &=& e^{\lambda{x^3}}\left.\exp\left(3x\lambda\frac{\partial^2}{\partial\xi^2}+\lambda\frac{\partial^3}{\partial\xi^3}\right)f(\xi+3x^2\lambda)\right|_{\xi=0}\\
   &=& e^{\lambda{x^3}}\exp\left(3x\lambda\frac{\partial^2}{\partial\zeta^2}+\lambda\frac{\partial^3}{\partial\zeta^3}\right)f(\zeta).
\end{eqnarray*}
Thus we obtain $\forall{f(\zeta)}\in{C^\infty}(\mathbb{R}^1)$,
\begin{equation}\label{y5.4}
e^{-\lambda{x^3}}f(\partial_x)e^{\lambda{x^3}}
=\exp\left(3x\lambda\frac{\partial^2}{\partial\zeta^2}+\lambda\frac{\partial^3}{\partial\zeta^3}\right)f(\zeta)\quad(\zeta=3x^2\lambda).
\end{equation}
Similarly,
\begin{equation}\label{y5.4+}
e^{-\lambda{x^2}}f(\partial_x)e^{\lambda{x^2}}=\exp\left(\lambda\frac{\partial^2}{\partial\zeta^2}\right)f(\zeta)\quad(\zeta=2x\lambda).
\end{equation}

In the special case when $f(\zeta)=\zeta^n$ for $n\in\mathbb{N}$, (\ref{y5.4}) can be written in the form:
\begin{equation}\label{y5.5}
e^{-\lambda{x^3}}\frac{d^n}{dx^n}\,e^{\lambda{x^3}}
=\sum^{[n/2]}_{k=0}\frac{\lambda^k}{k!}\left(3x+\frac{\partial}{\partial\zeta}\right)^k\frac{d^{2k}}{d\zeta^{2k}}\,\zeta^n\qquad(\zeta=3x^2\lambda).
\end{equation}
Letting $y=x^3$, then (\ref{y5.5}) gives the following operator equality:
\[\frac{d^n}{dx^n}\,e^{\lambda{x^3}}=\sum^{[n/2]}_{k=0}\sum^k_{j=0}\frac{3^{n-k-2j}n!}{j!\,(k-j)!\,(n-2k-j)!}x^{2n-3k-3j}\frac{d^{n-k-j}}{dy^{n-k-j}}\,e^{\lambda{y}}.\]
Thus we obtain (\ref{y5.3}) according to the Analytic Continuity Fundamental Theorem.
Similarly, we can obtain (\ref{y5}) from (\ref{y5.4+}).

\subsection{Algorithms of abstract operators}
\noindent

According to the Analytic Continuity Fundamental Theorem and Definition 2.4, by using the algebraic properties of two symbols $\cos(h\xi)$ and $\sin(h\xi)$,
Guang-Qing Bi has obtained the following three groups of algorithms of abstract operators just as the differential rules (\ref{3}) to (\ref{5}):

\textbf{Theorem BI1.} (See \cite[p. 9, Theorem 3, 4, and 6]{bi97}) Let $x\in\mathbb{R}^n,\;h\in\mathbb{R}^n$, $h\partial_x=\langle{h,\partial_x}\rangle$. Then for the abstract operators $\cos(h\partial_x)\,\mbox{and}\,\sin(h\partial_x)$ we have
\begin{equation}\label{y0}
  \cos(h\partial_x)f(x)=\Re[f(x+ih)],\quad\sin(h\partial_x)f(x)=\Im[f(x+ih)],
\end{equation}
$\forall{f(z)}\in{C}^\infty(\Omega),\,z=x+iy\in\Omega\subseteq\mathbb{C}^n$;

\parbox{11cm}{\begin{eqnarray*}\label{y1}
\sin(h\partial_x)(uv) &=& \cos(h\partial_x)v\cdot\sin(h\partial_x)u+\sin(h\partial_x)v\cdot\cos(h\partial_x)u,\\
\cos(h\partial_x)(uv) &=& \cos(h\partial_x)v\cdot\cos(h\partial_x)u-\sin(h\partial_x)v\cdot\sin(h\partial_x)u;
              \end{eqnarray*}}\hfill\parbox{2cm}{\begin{eqnarray}\end{eqnarray}}

\parbox{12cm}{\begin{eqnarray*}\label{y2}
                \sin(h\partial_x)\frac{u}{v} &=& \frac{\cos(h\partial_x)v\cdot\sin(h\partial_x)u-\sin(h\partial_x)v\cdot\cos(h\partial_x)u}
                {(\cos(h\partial_x)v)^2+(\sin(h\partial_x)v)^2}, \\
                \cos(h\partial_x)\frac{u}{v} &=& \frac{\cos(h\partial_x)v\cdot\cos(h\partial_x)u+\sin(h\partial_x)v\cdot\sin(h\partial_x)u}
                {(\cos(h\partial_x)v)^2+(\sin(h\partial_x)v)^2}.
              \end{eqnarray*}}\hfill\parbox{2cm}{\begin{eqnarray}\end{eqnarray}}

\textbf{Theorem BI2.} (See \cite[p. 9, Theorem 5]{bi97}) Let $h_0\in\mathbb{R},\;x(t)\in\mathbb{R}^n,\;t\in\mathbb{R}^1,\;X\in\mathbb{R}^n,\,Y\in\mathbb{R}^n$,
$Y\partial_X=\langle{Y,\partial_X}\rangle=Y_1\partial_{X_1}+\cdots+Y_n\partial_{X_n}$.
Then we have

\parbox{12cm}{\begin{eqnarray*}\label{y3}
 \sin(h_0\partial_t)f(x(t)) &=& \sin(Y\partial_X)f(X),\\
 \cos(h_0\partial_t)f(x(t)) &=& \cos(Y\partial_X)f(X),
\end{eqnarray*}}\hfill\parbox{2cm}{\begin{eqnarray}\end{eqnarray}}
where $X_j:=\cos(h_0\partial_t)x_j(t),\;Y_j:=\sin(h_0\partial_t)x_j(t),\;j=1,\cdots,n$.

In the special cases when $n=1$ and $n=2$, we have

\parbox{12cm}{\begin{eqnarray*}\label{yb1}
\sin\left(h_0\frac{d}{dt}\right)f(x(t)) &=& \sin\left(Y\frac{\partial}{\partial{X}}\right)f(X),\\
\cos\left(h_0\frac{d}{dt}\right)f(x(t)) &=& \cos\left(Y\frac{\partial}{\partial{X}}\right)f(X),
\end{eqnarray*}}\hfill\parbox{2cm}{\begin{eqnarray}\end{eqnarray}}
where
$Y:=\sin(h_0\partial_t)x(t),\;X:=\cos(h_0\partial_t)x(t),\;t\in\mathbb{R}^1,\;h_0\in\mathbb{R}$;

\parbox{12cm}{\begin{eqnarray*}\label{yb2}
\sin\left(h_0\frac{d}{dt}\right)f(x(t),y(t)) &=& \sin\left(Y_x\frac{\partial}{\partial{X}_x}+Y_y\frac{\partial}{\partial{X}_y}\right)f(X_x,X_y),\\
\cos\left(h_0\frac{d}{dt}\right)f(x(t),y(t)) &=& \cos\left(Y_x\frac{\partial}{\partial{X}_x}+Y_y\frac{\partial}{\partial{X}_y}\right)f(X_x,X_y),
\end{eqnarray*}}\hfill\parbox{2cm}{\begin{eqnarray}\end{eqnarray}}
where
$Y_x:=\sin(h_0\partial_t)x(t),\;X_x:=\cos(h_0\partial_t)x(t),\;Y_y:=\sin(h_0\partial_t)y(t),\;X_y:=\cos(h_0\partial_t)y(t)$.

\textbf{Theorem BI3.} (See \cite[p. 9, Theorem 7]{bi97}) Let $u=g(y)$ be a monotonic function on its domain. If $y=f(bx)$ is the inverse function of $bx=g(y)$ such that $g(f(bx))=bx$, where
$bx=b_1x_1+b_2x_2+\cdots+b_nx_n$ and $bh=b_1h_1+b_2h_2+\cdots+b_nh_n$, then $\sin(h\partial_x)f(bx)$ (denoted by $Y$) and
$\cos(h\partial_x)f(bx)$ (denoted by $X$) can be determined by the following set of equations:
\begin{equation}\label{y4}
    \left\{\begin{array}{l@{\qquad}l}\displaystyle
    \cos\left(Y\frac{\partial}{\partial{X}}\right)g(X)=bx,&x\in\mathbb{R}^n,\;b\in\mathbb{R}^n,\\\displaystyle
    \sin\left(Y\frac{\partial}{\partial{X}}\right)g(X)=bh, &h\in\mathbb{R}^n.
    \end{array}\right.
\end{equation}

\textbf{Proof.} According to the Analytic Continuity Fundamental Theorem and Definition 2.4, we have
\[e^{ih\partial_x}g(y)=g(f(bx+ibh))=g(e^{ih\partial_x}f(bx))=g(X+iY)=\exp\left(iY\frac{\partial}{\partial X}\right)g(X),\]
where
$y=f(bx),\;x\in\mathbb{R}^n,\;b\in\mathbb{R}^n;\;X:=\cos(h\partial_x)f(bx),\;Y:=\sin(h\partial_x)f(bx),\;h\in\mathbb{R}^n$.

According to Corollary 2.1, we have
\[e^{ih\partial_x}g(y)=\cos(h\partial_x)g(y)+i\sin(h\partial_x)g(y),\]
\[\exp\left(iY\frac{\partial}{\partial X}\right)g(X)=\cos\left(Y\frac{\partial}{\partial{X}}\right)g(X)+i\sin\left(Y\frac{\partial}{\partial{X}}\right)g(X).\]
Thus we obtain

\parbox{12cm}{\begin{eqnarray*}\label{ey3}
\sin(h\partial_x)g(y) &=& \sin\left(Y\frac{\partial}{\partial{X}}\right)g(X),\\
\cos(h\partial_x)g(y) &=& \cos\left(Y\frac{\partial}{\partial{X}}\right)g(X).
\end{eqnarray*}}\hfill\parbox{2cm}{\begin{eqnarray}\end{eqnarray}}
In the left side of (\ref{ey3}), if $u=g(y)$ is a monotonic function on its domain, then there exists function $y=f(bx)$ such that $g(f(bx))=bx$. Thus we have
$\sin(h\partial_x)g(y)=\sin(h\partial_x)bx=bh$ and $\cos(h\partial_x)g(y)=\cos(h\partial_x)bx=bx$ respectively. Theorem BI3 is proved.

Similarly, we can prove Theorem BI1 and Theorem BI2.

\textbf{Theorem 2.1.} Let $a=(a_1,\cdots,a_n)\in\mathbb{C}^n,\,a^{x\partial_x}:=a_1^{x_1\partial_1}a_2^{x_2\partial_2}\cdots{a}_n^{x_n\partial_n}$. Then
we have
\begin{equation}\label{15}
    a^{x\partial_x}f(x)=f(a_1x_1,\cdots,a_ix_i,\cdots,a_nx_n),\quad\forall{f}(x)\in{C}^\infty(\Omega),\;\Omega\in\mathbb{R}^n.
\end{equation}

\textbf{Proof.} By Definition 2.3, for the base function $x^\alpha,\,x\in\mathbb{R}^n,\,\alpha\in\mathbb{N}^n$ we have
\[a^{x\partial_x}x^\alpha=(a_1x_1)^{\alpha_1}(a_2x_2)^{\alpha_2}\cdots(a_nx_n)^{\alpha_n}.\]
Thus Theorem 2.1 is proved by using the Analytic Continuity Fundamental Theorem.

\textbf{Theorem 2.2.} Let $X:=(\rho_1\cos\theta_1,\cdots,\rho_n\cos\theta_n),\;Y:=(\rho_1\sin\theta_1,\cdots,\rho_n\sin\theta_n)$. Introducing the notation
\[\theta\,\rho\frac{\partial}{\partial\rho}:=\langle\theta,\rho\partial_\rho\rangle=\theta_1\,\rho_1\frac{\partial}{\partial\rho_1}+\cdots+
\theta_n\,\rho_n\frac{\partial}{\partial\rho_n},\]
\[Y\partial_X:=\langle{Y,\partial_X}\rangle=Y_1\frac{\partial}{\partial{X}_1}+\cdots+Y_n\frac{\partial}{\partial{X}_n}.\]
Then
$\forall{f}(\rho)\in{C}^\infty(\Omega),\;\rho\in\mathbb{R}^n,\;\theta\in\mathbb{R}^n$, we have

\parbox{12cm}{\begin{eqnarray*}\label{15*3}
                \cos\left(\theta\,\rho\frac{\partial}{\partial\rho}\right)f(\rho) &=& \cos(Y\partial_X)f(X), \\
                \sin\left(\theta\,\rho\frac{\partial}{\partial\rho}\right)f(\rho) &=& \sin(Y\partial_X)f(X),
\end{eqnarray*}}\hfill\parbox{2cm}{\begin{eqnarray}\end{eqnarray}}
which transforms the abstract operators taking $\rho\partial_\rho$ as the operator element into those taking $\partial_X$ as the operator element.

\textbf{Proof.} Let $a_1=e^{i\theta_1},\,\ldots,\,a_n=e^{i\theta_n}$. According to Theorem 2.1, we have
\[\exp\left(i\theta\,\rho\frac{\partial}{\partial\rho}\right)f(\rho)=f(\rho{e}^{i\theta})=f(X+iY)=\exp(iY\partial_X)f(X)\qquad (X\in\mathbb{R}^n,\;Y\in\mathbb{R}^n).\]
Considering
\[\cos\left(\theta\,\rho\frac{\partial}{\partial\rho}\right)f(\rho)+i\sin\left(\theta\,\rho\frac{\partial}{\partial\rho}\right)f(\rho)=
\cos(Y\partial_X)f(X)+i\sin(Y\partial_X)f(X),\]
thus we obtain (\ref{15*3}). Theorem 2.2 is proved.

\textbf{Corollary 2.4.} Let $f(z)\in{C^\infty}(\Omega),\,z=x+iy\in\Omega\subseteq\mathbb{C}^1$ be any analytic function. Then the harmonic functions $u(x,y)$ and $\upsilon(x,y)$ on complex plane can be expressed as

\parbox{12cm}{\begin{eqnarray*}\label{15**}
u(x,y)=\cos\left(y\frac{\partial}{\partial{x}}\right)f(x)\!\!
&=& \!\!\cos\left(\theta\rho\frac{\partial}{\partial\rho}\right)f(\rho),\\
\upsilon(x,y)=\sin\left(y\frac{\partial}{\partial{x}}\right)f(x)\!\!
&=& \!\!\sin\left(\theta\rho\frac{\partial}{\partial\rho}\right)f(\rho),
\end{eqnarray*}}\hfill\parbox{2cm}{\begin{eqnarray}\end{eqnarray}}
where $\rho=\sqrt{x^2+y^2},\,\theta=\arctan\frac{y}{x}$. Especially, if $y=kx$, then (\ref{15**}) gives that

\parbox{12cm}{\begin{eqnarray*}\label{15*}
\left.\cos\left(y\frac{\partial}{\partial{x}}\right)f(x)\right|_{y=kx}
 &=& \beta^{x\frac{\partial}{\partial{x}}}\cos\left(\alpha{x}\frac{\partial}{\partial{x}}\right)f(x),\\
\left.\sin\left(y\frac{\partial}{\partial{x}}\right)f(x)\right|_{y=kx}
 &=& \beta^{x\frac{\partial}{\partial{x}}}\sin\left(\alpha{x}\frac{\partial}{\partial{x}}\right)f(x),
\end{eqnarray*}}\hfill\parbox{2cm}{\begin{eqnarray}\end{eqnarray}}
where $\beta=\sqrt{1+k^2},\;\alpha=\arctan{k}$, which can be used to solve the boundary value problems of 2-dimensional Laplace equation on polygonal domains.

\subsection{The summation method of Fourier series}
\noindent

\textbf{Theorem 2.3.}  Let $f(x)\in{L^2}[-c,c]$ defined by the Fourier cosine series, and $g(x)\in{L^2}[-c,c]$ be that of the corresponding Fourier  sine series, namely
\[ f(x):=\sum^\infty_{n=0}a_n\cos\frac{n\pi{x}}{c}\quad\mbox{and}\quad g(x):=\sum^\infty_{n=0}a_n\sin\frac{n\pi{x}}{c},\]
where $c>0,\;x\in\Omega\subseteq\mathbb{R}^1$. If $S(t)$ is the sum function of the corresponding power series $\sum^\infty_{n=0}a_nt^n$, namely
\[S(t)=\sum^\infty_{n=0}a_nt^n\qquad(t\in\mathbb{R}^1,\;|t|<r,\;0<r<+\infty),\]
then we have the following trigonometric summation relationships:
\begin{equation}\label{12ct1}
 \left.\cos\left(\frac{\pi{x}}{c}\frac{\partial}{\partial{z}}\right)S(e^z)\right|_{z=0}=\sum^\infty_{n=0}a_n\cos\frac{n\pi{x}}{c};
\end{equation}
\begin{equation}\label{12ct2}
 \left.\sin\left(\frac{\pi{x}}{c}\frac{\partial}{\partial{z}}\right)S(e^z)\right|_{z=0}=\sum^\infty_{n=0}a_n\sin\frac{n\pi{x}}{c}.
\end{equation}
Here $\Omega$ can be uniquely determined by the detailed computation of the left-hand side of (\ref{12ct1}) and (\ref{12ct2}) respectively.

\textbf{Proof.} Theorem 2.3 can be proved easily by substituting $S(e^z)=\sum^\infty_{n=0}a_ne^{nz}$ into (\ref{12ct1}) and (\ref{12ct2}) respectively.

\textbf{Example 2.5.} Let $f(x)\in{L}^2([-c,c])$ be the square wave function defined by
\[f(x):=\left\{\begin{array}{r@{\qquad}l}
+1,&2mc-c/2<x<2mc+c/2,\\
-1,&2mc+c/2<x<2mc+3c/2,\end{array}\right.\]
where $m=0,\pm1,\pm2,\cdots$. Then $f(x)$ can be expressed in the form:
\[f(x)=\frac{4}{\pi}\left.\cos\left(\frac{\pi{x}}{c}\frac{\partial}{\partial{z}}\right)\arctan{e^z}\right|_{z=0}\in{C^\infty(\Omega)},\quad x\in\Omega:=\{x\in\mathbb{R}^1|\cos(\pi{x}/c)\neq0\}.\]

\textbf{Proof.} By making use of (\ref{y0}), we have

\parbox{12cm}{\begin{eqnarray*}\label{2a}
                \cos(h\partial_x)\sin{bx} &=& \cosh(bh)\sin{bx}, \\
                \sin(h\partial_x)\sin{bx} &=& \sinh(bh)\cos{bx}.
              \end{eqnarray*}}\hfill\parbox{2cm}{\begin{eqnarray}\end{eqnarray}}

\parbox{12cm}{\begin{eqnarray*}\label{2b}
                \cos(h\partial_x)\cos{bx} &=& \cosh(bh)\cos{bx}, \\
                \sin(h\partial_x)\cos{bx} &=& -\sinh(bh)\sin{bx}.
              \end{eqnarray*}}\hfill\parbox{2cm}{\begin{eqnarray}\end{eqnarray}}

Based on (\ref{2a}) and (\ref{2b}), by using (\ref{y2}) we obtain

\parbox{12cm}{\begin{eqnarray*}\label{2c}
                \cos(h\partial_x)\tan{bx} &=& \frac{\sin(2bx)}{\cosh(2bh)+\cos(2bx)}, \\
                \sin(h\partial_x)\tan{bx} &=& \frac{\sinh(2bh)}{\cosh(2bh)+\cos(2bx)}.
              \end{eqnarray*}}\hfill\parbox{2cm}{\begin{eqnarray}\end{eqnarray}}

By making use of (\ref{2c}) and (\ref{y4}), we have the following set of equations:
\begin{equation}\label{y4s}
    \left\{\begin{array}{l@{\qquad}l}\displaystyle
    \frac{\sin2X}{\cosh2Y+\cos2X} = bx,  &X=\cos(h\partial_x)\arctan{bx},\\\displaystyle
    \frac{\sinh2Y}{\cosh2Y+\cos2X} = bh, &Y=\sin(h\partial_x)\arctan{bx}.
    \end{array}\right.
\end{equation}

Solving $X$ and $Y$ from (\ref{y4s}), then
\[1+(bx)^2+(bh)^2=1+\frac{\sin^22X+\sinh^22Y}{(\cosh2Y+\cos2X)^2}=\frac{2\cosh2Y}{\cosh2Y+\cos2X}=\frac{2bh\cosh2Y}{\sinh2Y}\]
and
\[1-(bx)^2-(bh)^2=1-\frac{\sin^22X+\sinh^22Y}{(\cosh2Y+\cos2X)^2}=\frac{2\cos2X}{\cosh2Y+\cos2X}=\frac{2bx\cos2X}{\sin2X}.\]
So we have
\[\tanh2Y=\frac{2bh}{1+(bx)^2+(bh)^2}\quad\mbox{and}\quad\tan2X=\frac{2bx}{1-(bx)^2-(bh)^2}.\]
Thus we obtain

\parbox{12cm}{\begin{eqnarray*}\label{ly8}
                \sin(h\partial_x)\arctan{bx} &=& \frac{1}{2}\textrm{tanh}^{-1}\frac{2bh}{1+(bx)^2+(bh)^2}, \\
                \cos(h\partial_x)\arctan{bx} &=&
                \frac{1}{2}\arctan\frac{2bx}{1-(bx)^2-(bh)^2},
              \end{eqnarray*}}\hfill\parbox{2cm}{\begin{eqnarray}\end{eqnarray}}
where
$bx=b_1x_1+b_2x_2+\cdots+b_nx_n,\;bh=b_1h_1+b_2h_2+\cdots+b_nh_n$.

By making use of (\ref{yb1}) and (\ref{ly8}), we have
\begin{eqnarray*}
   & & \frac{4}{\pi}\left.\cos\left(\frac{\pi{x}}{c}\frac{\partial}{\partial{z}}\right)\arctan{e^z}\right|_{z=0}=\,
   \frac{4}{\pi}\left.\cos\left(Y\frac{\partial}{\partial{X}}\right)\arctan{X}\right|_{z=0}\\
   &=&\left.\frac{2}{\pi}\arctan\frac{2X}{1-(X^2+Y^2)}\right|_{z=0}\,=\,
   \frac{2}{\pi}\arctan\frac{2\cos(\pi{x}/c)}{1-\left(\cos^2(\pi{x}/c)+\sin^2(\pi{x}/c)\right)}\\
   &=&\left\{\begin{array}{r@{\qquad}l}
(2/\pi)\arctan(+\infty)=+1,&2mc-c/2<x<2mc+c/2\quad\;\;(\cos(\pi{x}/c)>0),\\
(2/\pi)\arctan(-\infty)=-1,&2mc+c/2<x<2mc+3c/2\quad(\cos(\pi{x}/c)<0),\end{array}\right.
\end{eqnarray*}
where $m=0,\pm1,\pm2,\cdots$. Thus Example 2.5 is proved.

More generally,  we have

\textbf{Theorem 2.4.} Let $x\in\Omega\subseteq\mathbb{R}^1$ be an open set. If $\forall{f(x)}\in{L}^2([-l,l])$ and $f(x+2l)=f(x)$ on $\Omega$, the analytic functions $S_+(t)$ and $S_-(t)$ are given by
\begin{equation}\label{f1'}
S_+(t):=\frac{1}{2l}\int^l_{-l}f(\xi)\frac{1-t^2}{1-2t\cos(\pi\xi/l)+t^2}d\xi
\end{equation}
and
\begin{equation}\label{f2'}
 S_-(t):=\frac{1}{l}\int^l_{-l}f(\xi)\frac{t\sin(\pi\xi/l)}{1-2t\cos(\pi\xi/l)+t^2}d\xi,
\end{equation}
respectively, then there exists $f_z(x)\in{C^\infty}(\mathbb{R}^1)$ with the following form
\begin{equation}\label{f3'}
f_z(x)=\cos\left(\frac{\pi{x}}{l}\frac{\partial}{\partial{z}}\right)S_+(e^z)+\sin\left(\frac{\pi{x}}{l}\frac{\partial}{\partial{z}}\right)S_-(e^z),\quad
-\infty<z<0
\end{equation}
such that
$$\left.f_z(x)\right|_{z=0}=\lim_{z\rightarrow0^-}f_z(x)\rightharpoonup{f(x)},\quad\forall{f(x)}\in{L}^2([-l,l]),$$
or, equivalently,
\begin{equation}\label{f2'd}
  f(x)=\left.\cos\left(\frac{\pi{x}}{l}\frac{\partial}{\partial{z}}\right)S_+(e^z)\right|_{z=0}
  +\left.\sin\left(\frac{\pi{x}}{l}\frac{\partial}{\partial{z}}\right)S_-(e^z)\right|_{z=0}\in{C^\infty(\Omega)}.
\end{equation}

\textbf{Proof.} By the algorithms (\ref{y2}), we have
\[\cos\left(\frac{\pi\xi}{l}\frac{\partial}{\partial{z}}\right)\frac{1+e^z}{1-e^z}=\frac{1-e^{2z}}{1-2e^{z}\cos(\pi\xi/l)+e^{2z}};\]
\[\sin\left(\frac{\pi\xi}{l}\frac{\partial}{\partial{z}}\right)\frac{e^z}{1-e^z}=\frac{e^z\sin(\pi\xi/l)}{1-2e^{z}\cos(\pi\xi/l)+e^{2z}}.\]
So $f_z(x)$ can be expressed as
\begin{eqnarray*}
   f_z(x) &=&
   \cos\left(\frac{\pi{x}}{l}\frac{\partial}{\partial{z}}\right)\left[\frac{1}{2l}\int^l_{-l}f(\xi)
   \cos\left(\frac{\pi\xi}{l}\frac{\partial}{\partial{z}}\right)\frac{1+e^z}{1-e^z}d\xi\right]\\
   & & +\,\sin\left(\frac{\pi{x}}{l}\frac{\partial}{\partial{z}}\right)\left[\frac{1}{l}\int^l_{-l}f(\xi)
   \sin\left(\frac{\pi\xi}{l}\frac{\partial}{\partial{z}}\right)\frac{e^z}{1-e^z}d\xi\right] \\
   &=& \frac{1}{2l}\int^l_{-l}f(\xi)\cos\left(\frac{\pi{x}}{l}\frac{\partial}{\partial{z}}\right)
   \cos\left(\frac{\pi\xi}{l}\frac{\partial}{\partial{z}}\right)\left[1+\frac{2e^z}{1-e^z}\right]d\xi\\
   & & +\,\frac{1}{l}\int^l_{-l}f(\xi)\sin\left(\frac{\pi{x}}{l}\frac{\partial}{\partial{z}}\right)
   \sin\left(\frac{\pi\xi}{l}\frac{\partial}{\partial{z}}\right)\frac{e^z}{1-e^z}d\xi = \frac{1}{2l}\int^l_{-l}f(\xi)d\xi\\
   & & +\,\frac{1}{l}\int^l_{-l}f(\xi)\left[\cos\left(\frac{\pi{x}}{l}\frac{\partial}{\partial{z}}\right)
   \cos\left(\frac{\pi\xi}{l}\frac{\partial}{\partial{z}}\right)+\sin\left(\frac{\pi{x}}{l}\frac{\partial}{\partial{z}}\right)
   \sin\left(\frac{\pi\xi}{l}\frac{\partial}{\partial{z}}\right)\right]\frac{e^z}{1-e^z}d\xi\\
   &=&
   \frac{1}{2l}\int^l_{-l}f(\xi)d\xi+\frac{1}{l}\int^l_{-l}f(\xi)\cos\left(\frac{\pi(x-\xi)}{l}\frac{\partial}{\partial{z}}\right)\frac{e^z}{1-e^z}d\xi\\
   &=& \frac{1}{2l}\int^l_{-l}f(\xi)\cos\left(\frac{\pi(x-\xi)}{l}\frac{\partial}{\partial{z}}\right)\frac{1+e^z}{1-e^z}d\xi\\
   &=& \frac{1}{2l}\int^l_{-l}f(\xi)\frac{1-e^{2z}}{1-2e^{z}\cos(\pi(x-\xi)/l)+e^{2z}}d\xi,\quad-\infty<z<0.
\end{eqnarray*}
\begin{eqnarray*}
 \lim_{z\rightarrow0^-}f_z(x) &=& \lim_{z\rightarrow0^-}\frac{1}{2l}\int^l_{-l}f(\xi)\frac{1-e^{2z}}{1-2e^{z}\cos(\pi(x-\xi)/l)+e^{2z}}d\xi\\
 &=&
\int^l_{-l}f(\xi)\lim_{z\rightarrow0^-}\frac{1}{2l}\frac{1-e^{2z}}{1-2e^{z}\cos(\pi(x-\xi)/l)+e^{2z}}d\xi
\rightharpoonup\int^l_{-l}f(\xi)\delta(x-\xi)d\xi\\
 &=& f(x).
\end{eqnarray*}
Thus Theorem 2.4 is proved.

\textbf{Remark 2.4.} In Theorem 2.4, $f(x)=\lim_{z\rightarrow0^-}f_z(x)\in{L}^2([-l,l])$, which may be discontinuous on real axis such as $f(x)$ is the square wave function in Example 2.5, but $f_z(x)$ corresponding to $f(x)$ is a $C^\infty$ function. In other words, if the mollifier $\rho_z(x)\in{C^\infty_0}(\mathbb{R}^1)$ is defined by
\[\rho_z(x):=\left\{\begin{array}{l@{\qquad}l}
\displaystyle\frac{1}{2l}\frac{1-e^{2z}}{1-2e^{z}\cos(\pi x/l)+e^{2z}},&|x|<l,\\
0,&|x|\geq l,\end{array}\right.\]
then $f_z(x)$ is given by the convolution relationship $f_z(x)=\rho_z(x)*f(x)$. This means that we can regularize $f(x)$ by the standard mollifier to get $f_z(x)$.
Therefore, for the abstract operators $g(\partial_x)$, if $\lim_{z\rightarrow0^-}g(\partial_x)f_z(x)\in{L}^2([-l,l])$, then $g(\partial_x)f(x)$ makes sense in a broad sense, which can be defined as
\begin{equation}\label{f4'}
 g(\partial_x)f(x):=\lim_{z\rightarrow0^-}g(\partial_x)f_z(x).
\end{equation}

\textbf{Theorem 2.5.}
 Let $t>0,\,\Re(s)>0$. If $\forall{f(x)}\in\mathscr{S}(\mathbb{R}^1)$, $F_+(s)$ and $F_-(s)$ are given by
\begin{equation}\label{f1a}
F_+(s):=\frac{1}{\pi}\mathcal{L}\int^\infty_{-\infty}f(\xi)\cos(t\xi)d\xi\quad\mbox{or}\quad F_+(s):=\frac{1}{\pi}\int^\infty_{-\infty}f(\xi)\frac{s}{\xi^2+s^2}d\xi
\end{equation}
and
\begin{equation}\label{f1b}
 F_-(s):=\frac{1}{\pi}\mathcal{L}\int^\infty_{-\infty}f(\xi)\sin(t\xi)d\xi\quad\mbox{or}\quad F_-(s):=\frac{1}{\pi}\int^\infty_{-\infty}f(\xi)\frac{\xi}{\xi^2+s^2}d\xi,
\end{equation}
respectively, where $\mathcal{L}$ is the Laplace transform, then we have
\begin{equation}\label{f3}
f(x)=\left.\cos\left(x\frac{\partial}{\partial{s}}\right)F_+(s)\right|_{s=0}+\left.\sin\left(-x\frac{\partial}{\partial{s}}\right)F_-(s)\right|_{s=0}.
\end{equation}

\textbf{Proof.} By (\ref{t4}) we have
\[\mathcal{L}\cos(t\xi)=\cos\left(\xi\frac{\partial}{\partial{s}}\right)\frac{1}{s},\quad
\mathcal{L}\sin(t\xi)=\sin\left(-\xi\frac{\partial}{\partial{s}}\right)\frac{1}{s}\qquad(\Re(s)>0).\]
So by using (\ref{f1a}) and (\ref{f1b}), then
$\forall{f(x)}\in\mathscr{S}(\mathbb{R}^1)$, there exists function $f_s(x)\in{C^\infty}(\mathbb{R}^1)$, which can be expressed as
\begin{eqnarray*}
   f_s(x) &=& \cos\left(x\frac{\partial}{\partial{s}}\right)F_+(s)+\sin\left(-x\frac{\partial}{\partial{s}}\right)F_-(s) \\
   &=& \cos\left(x\frac{\partial}{\partial{s}}\right)\left[\frac{1}{\pi}\int^\infty_{-\infty}f(\xi)\mathcal{L}\cos(t\xi)d\xi\right]
   +\sin\left(-x\frac{\partial}{\partial{s}}\right)\left[\frac{1}{\pi}\int^\infty_{-\infty}f(\xi)\mathcal{L}\sin(t\xi)d\xi\right]\\
   &=& \frac{1}{\pi}\int^\infty_{-\infty}f(\xi)\left[\cos\left(x\frac{\partial}{\partial{s}}\right)
   \cos\left(\xi\frac{\partial}{\partial{s}}\right)+\sin\left(-x\frac{\partial}{\partial{s}}\right)
   \sin\left(-\xi\frac{\partial}{\partial{s}}\right)\right]\frac{1}{s}\,d\xi\\
   &=& \frac{1}{\pi}\int^\infty_{-\infty}f(\xi)\cos\left((x-\xi)\frac{\partial}{\partial{s}}\right)\frac{1}{s}\,d\xi
   =\frac{1}{\pi}\int^\infty_{-\infty}f(\xi)\frac{s}{(x-\xi)^2+s^2}d\xi,\;\Re(s)>0.
\end{eqnarray*}
Therefore, we have
\begin{eqnarray*}
 \lim_{s\rightarrow0^+}f_s(x) &=& \lim_{s\rightarrow0^+}\frac{1}{\pi}\int^\infty_{-\infty}f(\xi)\frac{s}{(x-\xi)^2+s^2}d\xi\\
 &=& \int^\infty_{-\infty}f(\xi)\lim_{s\rightarrow0^+}\frac{1}{\pi}\frac{s}{(x-\xi)^2+s^2}d\xi\rightharpoonup\int^\infty_{-\infty}f(\xi)\delta(x-\xi)d\xi=f(x).
\end{eqnarray*}
Thus Theorem 2.5 is proved.

Substituting (\ref{f3}) into (\ref{f1a}) and (\ref{f1b}) respectively, then Theorem 2.5 gives that

\textbf{Corollary 2.5.}  Let $t>0,\,F(s):=\mathcal{L}f(t)$. Then $\forall{f(t)}\in\mathscr{S}(\mathbb{R}^1)$, the inverse of Laplace transform can be determined by
\begin{equation}\label{f1+}
\mathcal{L}^{-1}F(s)=\frac{2}{\pi}\int^\infty_0 \!\left.\cos\left(\xi\frac{\partial}{\partial{s}}\right)F(s)\right|_{s=0}\cos(t\xi)d\xi
\end{equation}
and
\begin{equation}\label{f1++}
\mathcal{L}^{-1}F(s)=\frac{2}{\pi}\int^\infty_0\!\left.\sin\left(-\xi\frac{\partial}{\partial{s}}\right)F(s)\right|_{s=0}\sin(t\xi)d\xi,
\end{equation}
respectively.

In fact, if $f(x)$ is the rational proper functions, then the integral representations (\ref{f1a}) and (\ref{f1b}) are absolutely convergent. Therefore, by Theorem 2.5 we have

\textbf{Corollary 2.6.} Let $f(x),\,x\in\mathbb{R}^1$ be the rational proper functions. If analytic functions $F_+(s)$ satisfies the following operator equation
\begin{equation}\label{f3+}
\left.\cos\left(x\frac{\partial}{\partial{s}}\right)F_+(s)\right|_{s=0}=\frac{1}{2}(f(x)+f(-x)),
\end{equation}
then we have
\begin{equation}\label{f1+3}
\int^\infty_{-\infty}f(x)\cos(tx)dx=\pi\mathcal{L}^{-1}F_+(s),\quad t>0.
\end{equation}
If analytic functions $F_-(s)$ satisfies the following operator equation
\begin{equation}\label{f3++}
\left.\sin\left(-x\frac{\partial}{\partial{s}}\right)F_-(s)\right|_{s=0}=\frac{1}{2}(f(x)-f(-x)),
\end{equation}
then we have
\begin{equation}\label{f1+3+}
\int^\infty_{-\infty}f(x)\sin(tx)dx=\pi\mathcal{L}^{-1}F_-(s),\quad t>0.
\end{equation}

\subsection{Integral representations for abstract operators}
\noindent

\textbf{Theorem BI4.} Let $P(\partial_x)$ be a $m$-order constant coefficient linear partial differential operators, $t\in\mathbb{R}^1$ with $t>0$. If there exist
$a_1,\,a_2,\,\cdots,\,a_k$ of real and partial differential operators $A_1,\,A_2,\,\cdots,\,A_k$ of the order less than $[(m+1)/2]$ such that
$P(\partial_x)\equiv{a}_1A_1^2+a_2A_2^2+\cdots+a_kA_k^2$ for $k=2\nu+3,\;\nu=0,\,1,\,2,\cdots$, then
$\forall{f(x)}\in{C^\infty}(\Omega),\,\Omega\in\mathbb{R}^n$, we have the following operator relationships, which was given by Guang-Qing Bi \cite[p. 12, Theorem 14]{bi97}:

\begin{eqnarray}\label{y6}
\lefteqn{\frac{\sinh\left(tP(\partial_x)^{1/2}\right)}{P(\partial_x)^{1/2}}f(x)=}\nonumber\\
   & &t\underbrace{\int^t_0tdt\cdots}_\nu\int^t_0\,tdt\frac{(P(\partial_x))^{\nu}}{2^{\nu+2}\pi^{\nu+1}} \int^\pi_{-\pi}\underbrace{\int^\pi_0\cdots}_{k-2}\int^\pi_0e^{\eta_1a_1^{1/2}A_1+\cdots+\eta_ka_k^{1/2}A_k}f(x)\,d\sigma_k\nonumber\\
   & & {}+\sum^{\nu-1}_{i=0}\frac{t^{2i+1}}{(2i+1)!}(P(\partial_x))^{i}f(x).
\end{eqnarray}
Similarly,
\begin{eqnarray}\label{y6'}
\lefteqn{\frac{\sin\left({t}P(\partial_x)^{1/2}\right)}{P(\partial_x)^{1/2}}f(x)=}\nonumber\\
   & & t\underbrace{\int^t_0tdt\cdots}_\nu\int^t_0tdt\frac{2^{\nu+1}(-P(\partial_x))^{\nu}}{\pi^{\nu+1}}
  \underbrace{\int^{\pi/2}_0\!\!\cdots}_{k-1}\int^{\pi/2}_0\!\!\cos(\eta_1a_1^{1/2}A_1)\cdots\cos(\eta_ka_k^{1/2}A_k)f(x)\,d\sigma_k\nonumber\\
   & & {}+\sum^{\nu-1}_{i=0}\frac{t^{2i+1}}{(2i+1)!}(-P(\partial_x))^{i}f(x).
\end{eqnarray}
Here $\eta\in\mathbb{R}^k$ is the integral variable and
\begin{eqnarray*}
  \eta_1 &=& t\cos\theta_1,\\
  \eta_2 &=& t\sin\theta_1\cos\theta_2,\\
  \eta_3 &=& t\sin\theta_1\sin\theta_2\cos\theta_3,\\
         &\cdots&\\
  \eta_p &=& t\sin\theta_1\sin\theta_2\cdots\sin\theta_{p-1}\cos\theta_p,\\
  \eta_{p+1} &=& t\sin\theta_1\sin\theta_2\cdots\sin\theta_p\cos\phi,\\
  \eta_{p+2} &=&
  \eta_k\;=\;t\sin\theta_1\sin\theta_2\cdots\sin\theta_p\sin\phi;
\end{eqnarray*}
\[d\sigma_k=\sin^{k-2}\theta_1\sin^{k-3}\theta_2\cdots\sin\theta_{k-2}d\theta_1d\theta_2\cdots{d}\theta_{k-2}d\phi.\]

\textbf{Proof.} In (\ref{y6}), let $f(x)=e^{\xi{x}},\,x\in\mathbb{R}^n,\xi\in\mathbb{R}^n$, and the symbols of the partial differential operators $A_j,j=1,2,\cdots,k$ be
denoted by $\chi_j(\xi),\;\beta_j:=a_j^{1/2}\chi_j(\xi)$. Then
(\ref{y6}) degenerates to its characteristic equation by Definition 2.4, namely
\begin{eqnarray}\label{y6z}
\frac{\sinh\left(tP(\xi)^{1/2}\right)}{P(\xi)^{1/2}}&=&
t\underbrace{\int^t_0tdt\cdots}_\nu\int^t_0\,tdt\frac{(P(\xi))^{\nu}}{2^{\nu+2}\pi^{\nu+1}}
\int^\pi_{-\pi}\underbrace{\int^\pi_0\cdots}_{k-2}\int^\pi_0e^{\eta_1\beta_1+\cdots+\eta_k\beta_k}d\sigma_k\nonumber\\
   & & +\,\sum^{\nu-1}_{i=0}\frac{t^{2i+1}}{(2i+1)!}(P(\xi))^i\quad(P(\xi)=\beta^2_1+\beta^2_2+\cdots+\beta^2_k).
\end{eqnarray}
According to the Analytic Continuity Fundamental Theorem, we only need to prove (\ref{y6z}). Solving the integral on a hypersphere on
the right side of (\ref{y6z}), we have
\begin{eqnarray*}
\frac{\sinh\left(tP(\xi)^{1/2}\right)}{P(\xi)^{1/2}}=t\underbrace{\int^t_0tdt\cdots}_\nu\int^t_0tdt
\sum^\infty_{j=0}\frac{(P(\xi))^{\nu+j}t^{2j}}{(2j)!!(2j+2\nu+1)!!}+\sum^{\nu-1}_{i=0}\frac{t^{2i+1}(P(\xi))^i}{(2i+1)!}.
\end{eqnarray*}
Then it is proved by the termwise integration of the infinite series on the right side of the equality. Similarly, we have (\ref{y6'}). Thus Theorem BI4 is proved.

When $A_1,\,\cdots,\,A_k$ in the right-hand of (\ref{y6}) and (\ref{y6'}) are one order partial differential operators, then the abstract operators
$e^{\eta_1a_1^{1/2}A_1+\cdots+\eta_ka_k^{1/2}A_k}$ and $\cos(\eta_1a_1^{1/2}A_1)\cdots\cos(\eta_ka_k^{1/2}A_k)$ are one of the following five simplest operators:
\[\exp(h\partial_x),\;\sin(h\partial_x),\;\cos(h\partial_x),\;\sinh(h\partial_x)\;\mbox{and}\;\cosh(h\partial_x).\]

It is easily seen from (\ref{y6}) of Theorem BI4 that

\textbf{Corollary 2.7.} Let $\Delta:=\sum^n_{k=1}\partial^2_{x_k}$ be the n-dimensional Laplacian, $a>0$ be the real parameter. If $n-2=2\nu+1,\,\nu\in\mathbb{N}_0:=\{0,1,2,\ldots\}$, then we have
\begin{eqnarray}\label{37}
\frac{\sinh\left(at\Delta^{1/2}\right)}{a\Delta^{1/2}}f(x)&=&
t\underbrace{\int^t_0tdt\cdots}_\nu\int^t_0\frac{(a^{2}\Delta)^{\nu}}{S_n}\int_{S_n}f(\xi)\,dS_n\,tdt\nonumber\\
   & & +\,\sum^{\nu-1}_{i=0}\frac{t^{2i+1}}{(2i+1)!}(a^{2}\Delta)^{i}f(x),\quad\forall{f(x)}\in{C}^{2\nu}(\Omega).
\end{eqnarray}
Here $t\in\mathbb{R}^1$ with $t>0$, $S_n=2(2\pi)^{\nu+1}(at)^{n-1}$. $\xi\in\mathbb{R}^n$ is the integral variable on the hypersphere
$(\xi_1-x_1)^2+(\xi_2-x_2)^2+\cdots+(\xi_n-x_n)^2=(at)^2$, and $dS_n$ is its surface element.

If $A_1,\,\cdots,\,A_k$ in (\ref{y6}) and (\ref{y6'}) are partial differential operators of the order great than 1, then the order can
be lowered by taking following result:

\textbf{Theorem BI5.} (See \cite[p. 11, Theorem 13]{bi97}) Let $x\in\Omega\subseteq\mathbb{R}^n$, $P(\partial_x)$ be the constant coefficient partial differential operators of any order. If
$f(x)\in{C}^\infty(\Omega)$ is a function for which the integral in the follow formula is finite, then
\begin{equation}\label{y7}
    e^{\lambda
    P(\partial_x)}f(x)=\frac{1}{2\sqrt{\pi}}\int^\infty_{-\infty}e^{-\zeta^2/4}e^{\lambda^{1/2}\zeta
    P(\partial_x)^{1/2}}f(x)d\zeta,\quad\lambda\in\mathbb{C}.
\end{equation}

By using Corollary 2.1, we deduce from (\ref{y7}) easily that

\textbf{Example 2.6.} Let $h_{\lambda,a}(\zeta)=a\lambda+\sqrt{\lambda/2}\zeta,\;\forall\lambda,a\in\mathbb{C},\;x\in\mathbb{R}^n,\;t\in\mathbb{R}^1$.
Then we have the following operator relationships:
\[\exp\left(-a^2t\frac{\partial^2}{\partial{x}_j^2}\right)g(x)=
\frac{1}{2\sqrt{\pi}}\int^\infty_{-\infty}e^{-\frac{\zeta^2}{4}}\cos\left(a\sqrt{t}\,\zeta\frac{\partial}{\partial{x_j}}\right)g(x)d\zeta;\]
\[\cos\left(\lambda\frac{\partial^2}{\partial{x}_j^2}\right)g(x)=\frac{1}{2\sqrt{\pi}}\int^\infty_{-\infty}e^{-\frac{\zeta^2}{4}}
e^{\sqrt{\lambda/2}\,\zeta\frac{\partial}{\partial{x}_j}}\cos\left(\sqrt{\frac{\lambda}{2}}\,\zeta\frac{\partial}{\partial{x}_j}\right)g(x)d\zeta;\]
\[\sin\left(\lambda\frac{\partial^2}{\partial{x}_j^2}\right)g(x)=\frac{1}{2\sqrt{\pi}}\int^\infty_{-\infty}e^{-\frac{\zeta^2}{4}}
e^{\sqrt{\lambda/2}\,\zeta\frac{\partial}{\partial{x}_j}}\sin\left(\sqrt{\frac{\lambda}{2}}\,\zeta\frac{\partial}{\partial{x}_j}\right)g(x)d\zeta;\]
\[\sin\left(a\lambda\frac{\partial}{\partial{x}_j}+\lambda\frac{\partial^2}{\partial{x}_j^2}\right)g(x)
=\frac{1}{2\sqrt{\pi}}\int^\infty_{-\infty}e^{-\frac{\zeta^2}{4}}
e^{\sqrt{\lambda/2}\,\zeta\frac{\partial}{\partial{x}_j}}\sin\left(h_{\lambda,a}(\zeta)\frac{\partial}{\partial{x}_j}\right)g(x)d\zeta.\]

\section{Main results}
\noindent\setcounter{equation}{0}
\renewcommand\theequation{3.\arabic{equation}}

Solving the ordinary or partial differential equations, is constructing the algorithms of the inverse operators of ordinary or partial differential operators. In this section we will see how abstract operators can be applied in solving partial differential equations without the experience of Fourier transform and its inversion process.

\subsection{The Laplace transform method for solving initial value problem of $n+1$-dimensional partial differential equations}
\noindent

In terms of abstract operators, solving the initial value problem of $n+1$ dimensional partial differential equations is similar to solving the ordinary differential equations with respect to the variable $t$, thus we can introduce the Laplace transform to further simplify the solving process. For this reason, we first need to use abstract operators to extend the mathematical concepts of Laplace transform.

Let $s\in\mathbb{C}^1$, $\lambda$ be the complex parameter. If taking $a,b\in\mathbb{R}$ as the characters of the base functions $e^{as}$ and $e^{bs}$ respectively, then we can combine $e^{as}$ and $e^{bs}$ with the following characteristic equation from Taylor formula (\ref{2})
\[f(\lambda a+\lambda b)=
\sum_{k=0}^\infty\frac{(\lambda a)^k}{k!}f^{(k)}(\lambda b),\quad|\lambda|<R\]
to make the following operator equality for the base functions $e^{as}$ and $e^{bs}$:
\[f\!\left(\lambda\frac{\partial}{\partial{s}}\right)(e^{as}e^{bs})=
\sum_{k=0}^\infty\frac{\lambda^k}{k!}\,\frac{\partial^ke^{as}}{\partial{s^k}}\,f^{(k)}\!\left(\lambda\frac{\partial}{\partial{s}}\right)e^{bs},\quad|\lambda|<R.\]
According to the Analytic Continuity Fundamental Theorem, we can obtain that

\textbf{Theorem 3.1.} Let $s\in\Omega\subset\mathbb{C}^1$. Suppose that $\lambda$ is the complex parameter. For the abstract operators $f^{(k)}(\lambda\partial_s),\,k=0,1,2,\cdots$, if there are two analytic functions $v(s),u(s)\in{C^\infty}(\Omega)$ such that the infinite series on the right side of (\ref{t2}) uniform convergent for $|\lambda|<R$, then it will uniform converges to the left side of this equality, namely
\begin{equation}\label{t2}
f\!\left(\lambda\frac{\partial}{\partial{s}}\right)(vu)=
\sum_{k=0}^\infty\frac{\lambda^k}{k!}\,\frac{\partial^kv}{\partial{s^k}}\,f^{(k)}\!\left(\lambda\frac{\partial}{\partial{s}}\right)u,\quad|\lambda|<R.
\end{equation}

When $\lambda=-1,\;v=s,\;u=1/s$ for $\Re(s)>0$, (\ref{t2}) gives that
\begin{equation}\label{t3}
f'\left(-\frac{\partial}{\partial{s}}\right)\frac{1}{s}=sf\left(-\frac{\partial}{\partial{s}}\right)\frac{1}{s}-f(0),
\end{equation}
which is equivalent to $\mathcal{L}f'(t)=s\mathcal{L}f(t)-f(0)$. The symbol $\mathcal{L}$ is the Laplace transform, which acts on functions $f(t)$ and generates a new function $F(s)=\mathcal{L}f(t)$. Thus we have

\textbf{Corollary 3.1.} Suppose that $f(t)$ is a real or complex valued function of the (time) variable $t>0$ and $s$ is a real or complex parameter. We can also define the Laplace transform of $f(t)$ as
\begin{equation}\label{t4}
F(s)=\mathcal{L}f(t):=f\left(-\frac{\partial}{\partial{s}}\right)\frac{1}{s},\quad\;\Re(s)>0.
\end{equation}

It is easily seen from (\ref{t4}) that
\begin{equation}\label{t5}
\mathcal{L}[g(t)f(t)]=g\left(-\frac{\partial}{\partial{s}}\right)F(s)\quad(F(s)=\mathcal{L}f(t));
\end{equation}
\begin{equation}\label{t4d}
f(t)=\mathcal{L}^{-1}F(s)=\mathcal{L}^{-1}f\left(-\frac{\partial}{\partial{s}}\right)\frac{1}{s},\quad\;\Re(s)>0,\;t\in\mathbb{R}_+^1:=\{t\in\mathbb{R}^1|t>0\}.
\end{equation}
Thus the symbols of the abstract operators on the right side of (\ref{t4}) and (\ref{t5}) can be further described by using conditions of the Laplace transform.

\textbf{Theorem 3.2.} Let $m\geq1$, $P(\partial_x)$ be an any order partial differential equations. Then we have
\begin{equation}\label{32}
    \left\{\begin{array}{l@{\qquad}l}\displaystyle
    \left(\frac{\partial^2}{\partial{t^2}}-P(\partial_x)\right)^mu=f(x,t),&x\in\mathbb{R}^n,\;t\in\mathbb{R}_+^1,\\\displaystyle
    \left.\frac{\partial^ru}{\partial{t^r}}\right|_{t=0}=\varphi_r(x),&r=0,1,2,\ldots,2m-1.
    \end{array}\right.
\end{equation}
\begin{eqnarray}\label{32'}
  u(x,t) &=&
  \int^t_0\int^{t-\tau}_0\frac{\left((t-\tau)^2-\tau'^2\right)^{m-2}}{(2m-2)!!\,(2m-4)!!}\,
  \frac{\sinh\left(\tau'P(\partial_x)^{1/2}\right)}{P(\partial_x)^{1/2}}\,f(x,\tau)\,\tau'd\tau'\,d\tau \nonumber\\
   & & +\,\sum^{m-1}_{k=0}(-1)^k{m\choose{k}}P(\partial_x)^k\sum^{2m-1-2k}_{r=0}\frac{\partial^{2m-1-2k-r}}{\partial t^{2m-1-2k-r}}\int^t_0
   \frac{(t^2-\tau^2)^{m-2}\tau}{(2m-2)!!\,(2m-4)!!}\nonumber\\
   & & \times\,\frac{\sinh\left(\tau P(\partial_x)^{1/2}\right)}{P(\partial_x)^{1/2}}\,\varphi_r(x)\,d\tau,
\end{eqnarray}
which was given by Guang-Qing Bi and Yue-Kai Bi \cite[p. 513, Theorem 1]{bi11}.

\textbf{Proof.} Taking the Laplace transform of both sides of the partial differential equations and considering the initial condition gives
\begin{equation}\label{32'bh}
  \sum^m_{k=0}(-1)^k{m\choose{k}}P(\partial_x)^k\left(s^{2m-2k}U(x,s)-\sum^{2m-1-2k}_{r=0}s^{2m-1-2k-r}\varphi_r(x)\right)=F(x,s),
\end{equation}
where $U(x,s):=\mathcal{L}u(x,t),\;F(x,s):=\mathcal{L}f(x,t)$.

We need to introduce the abstract operators $G_m(P(\partial_x),t)$ defined by
$$G_m(P(\partial_x),t):=\mathcal{L}^{-1}\frac{1}{(s^2-P(\partial_x))^m}.$$

Suppose that $f(t)$ is a real or complex valued function of the (time) variable $t>0$ and $s$ is a real or complex parameter. Then we have
\[\mathcal{L}\,\left(\int^t_0\!\cdot\,{t}dt\right)^{m-1}\!\!f(t)=\left(-\frac{1}{s}\frac{\partial}{\partial{s}}\right)^{m-1}\!\!\mathcal{L}f(t),\quad{m}\geq1.\]
Letting $f(t)=\sin{bt},\,b\in\mathbb{C},\;t\in\mathbb{R}_+$, we have
\[\left(\int^t_0\cdot\,tdt\right)^{m-1}\!\!\sin{bt}=\mathcal{L}^{-1}\left(-\frac{1}{s}\frac{\partial}{\partial{s}}\right)^{m-1}\!\!\frac{b}{s^2+b^2}
=\mathcal{L}^{-1}\frac{2^{m-1}(m-1)!}{(s^2+b^2)^m}b.\]
Letting $b=iP(\xi)^{1/2},\;\xi\in\mathbb{R}^n$, we have
\[\mathcal{L}^{-1}\frac{1}{(s^2-P(\xi))^m}=\frac{1}{(2m-2)!!}\left(\int^t_0\cdot\,tdt\right)^{m-1}\frac{\sinh\left(tP(\xi)^{1/2}\right)}{P(\xi)^{1/2}}.\]
Taking this one as the characteristic equation, according to the Analytic Continuity Fundamental Theorem, we have
\begin{equation}\label{33}
G_m(P(\partial_x),t)=\mathcal{L}^{-1}\frac{1}{(s^2-P(\partial_x))^m}=
\frac{1}{(2m-2)!!}\left(\int^t_0\cdot\,tdt\right)^{m-1}\frac{\sinh\left(tP(\partial_x)^{1/2}\right)}{P(\partial_x)^{1/2}}.
\end{equation}
It is easily seen from (\ref{33}) that
$$\left.\frac{\partial^k}{\partial t^k}G_m(P(\partial_x),t)\right|_{t=0}=0\qquad(k=0,1,2,\ldots,2m-2).$$
We can easily derive the following integral relationship
\begin{equation}\label{34}
\left(\int^t_a\cdot\,tdt\right)^{m}f(t)=\underbrace{\int^t_atdt\cdots}_{m}\int^t_af(t)\,tdt=\int^t_a\frac{(t^2-\tau^2)^{m-1}}{(2m-2)!!}f(\tau)\,\tau d\tau.
\end{equation}
Applying (\ref{34}) to (\ref{33}), we have the expression of abstract operators $G_m(P(\partial_x),t)$:
\begin{equation}\label{35}
G_m(P(\partial_x),t)g(x)=
\int^t_0\frac{(t^2-\tau'^2)^{m-2}}{(2m-2)!!\,(2m-4)!!}\,\frac{\sinh\left(\tau'P(\partial_x)^{1/2}\right)}{P(\partial_x)^{1/2}}g(x)\,\tau'd\tau'.
\end{equation}

Solving $U(x,s)$ from (\ref{32'bh}), its inverse transform is
\begin{eqnarray}\label{32'bi}
 u(x,t)&=& \mathcal{L}^{-1}U(x,s)=\mathcal{L}^{-1}\frac{1}{(s^2-P(\partial_x))^m}F(x,s)\nonumber\\
  & &+\,\mathcal{L}^{-1}\sum^{m-1}_{k=0}(-1)^k{m\choose{k}}P(\partial_x)^k\sum^{2m-1-2k}_{r=0}\frac{s^{2m-1-2k-r}}{(s^2-P(\partial_x))^m}\varphi_r(x)\nonumber\\
  &=& G_m(P(\partial_x),t)\ast f(x,t)\nonumber\\
  & &+\,\sum^{m-1}_{k=0}(-1)^k{m\choose{k}}P(\partial_x)^k\sum^{2m-1-2k}_{r=0}\frac{\partial^{2m-1-2k-r}}{\partial t^{2m-1-2k-r}}G_m(P(\partial_x),t)\varphi_r(x).
\end{eqnarray}
Applying (\ref{35}) to (\ref{32'bi}), thus Theorem 3.2 is proved.

It is extremely complex to prove Theorem 3.2 even if $\varphi_j(x)=0,\,j=0,1,\ldots,2m-1$ without using the Laplace transform (See Guang-Qing Bi \cite[p. 89, Theorem 1]{bi01}).

In 1999, Guang-Qing Bi \cite[p. 86, Theorem 2 and 3]{bi99} has obtained the following results:

\textbf{Theorem BI5.} Let $a_1,a_2,\ldots,a_m$ be any real or complex numbers different from each other, $P(\partial_x)$ be a partial differential operators of any order. Then we have
\begin{equation}\label{I1}
    \left\{\begin{array}{l@{\qquad}l}\displaystyle
    \prod^m_{i=1}(\frac{\partial}{\partial{t}}-a_iP(\partial_x))u=f(x,t),&x\in\mathbb{R}^n,\;t\in\mathbb{R}_+^1,\;m\geq1,\\\displaystyle
    \left.\frac{\partial^ju}{\partial{t^j}}\right|_{t=0}=0,&j=0,1,2,\ldots,m-1.
    \end{array}\right.
\end{equation}
\begin{equation}\label{I1'}
u(x,t)=\int^t_0\int^{t-\tau}_0\frac{(t-\tau-\tau')^{m-2}}{(m-2)!}\sum^m_{j=1}\frac{a_j^{m-1}}{\prod^m_{i=1\atop
i\neq{j}}(a_j-a_i)} e^{\tau'a_jP(\partial_x)}f(x,\tau)\,d\tau'd\tau;
\end{equation}

\textbf{Theorem BI6.} Let $a_1,a_2,\ldots,a_m$ be any positive real numbers different from each other, $P(\partial_x)$ be a partial differential operators of any order. Then we have
\begin{equation}\label{I2}
    \left\{\begin{array}{l@{\qquad}l}\displaystyle
    \prod^m_{i=1}(\frac{\partial^2}{\partial{t^2}}-a_i^2P(\partial_x))u=f(x,t),&x\in\mathbb{R}^n,\;t\in\mathbb{R}_+^1,\;m\geq1,\\\displaystyle
    \left.\frac{\partial^ju}{\partial{t^j}}\right|_{t=0}=0,&j=0,1,2,\ldots,2m-1.
    \end{array}\right.
\end{equation}
\begin{equation}\label{I2'}
u(x,t)=\int^t_0\int^{t-\tau}_0\frac{(t-\tau-\tau')^{2m-3}}{(2m-3)!}\sum^m_{j=1}\frac{a_j^{2m-2}}{\prod^m_{i=1
\atop i\neq{j}}(a_j^2-a_i^2)}
\frac{\sinh(\tau'a_jP(\partial_x)^{1/2})}{a_jP(\partial_x)^{1/2}}f(x,\tau)\,d\tau'd\tau.
\end{equation}

On this basis, by using the abstract operators and Laplace transform we have obtained the following theorems:

\textbf{Theorem 3.3.} Let $a_1,a_2,\ldots,a_m$ be the real or complex roots different from each other for any algebraic equation
$b_0+b_1\chi+b_2\chi^2+\cdots+b_m\chi^m=0$, and $P(\partial_x,\partial_t)$ be the partial differential operators defined by
$$P(\partial_x,\partial_t):=\sum^m_{k=0}b_kP(\partial_x)^{m-k}\frac{\partial^k}{\partial{t^k}}=\prod^m_{i=1}\left(\frac{\partial}{\partial t}-a_iP(\partial_x)\right),
\quad{x}\in\mathbb{R}^n,\;t\in\mathbb{R}_+^1,\;m\geq1.$$
Here $P(\partial_x)$ is a partial differential operators of any order. Then we have
\begin{equation}\label{j1}
    \left\{\begin{array}{l@{\qquad}l}\displaystyle
    P(\partial_x,\partial_t)u=f(x,t),&x\in\mathbb{R}^n,\;t\in\mathbb{R}_+^1,\\\displaystyle
    \left.\frac{\partial^ru}{\partial{t^r}}\right|_{t=0}=\varphi_r(x),&r=0,1,2,\ldots,m-1.
    \end{array}\right.
\end{equation}
\begin{eqnarray}\label{j1'}
  u(x,t) &=&
  \int^t_0\int^{t-\tau}_0\frac{(t-\tau-\tau')^{m-2}}{(m-2)!}\sum^m_{j=1}\frac{a_j^{m-1}}{\prod^m_{i=1 \atop i\neq{j}}(a_j-a_i)}\,
e^{\tau'a_jP(\partial_x)}f(x,\tau)\,d\tau'd\tau\nonumber\\
   & & +\,\sum^m_{k=0}b_kP(\partial_x)^{m-k}\sum^{k-1}_{r=0}\frac{\partial^{k-1-r}}{\partial t^{k-1-r}}\int^t_0
   \frac{(t-\tau)^{m-2}}{(m-2)!}\nonumber\\
   & & \times\,\sum^m_{j=1}\frac{a_j^{m-1}}{\prod^m_{i=1 \atop i\neq{j}}(a_j-a_i)}\,e^{\tau{a_j}P(\partial_x)}\,\varphi_r(x)\,d\tau.
\end{eqnarray}

\textbf{Theorem 3.4.} Let $a_1,a_2,\ldots,a_m$ be any positive real numbers different from each other such that
$\sum^m_{k=0}b_{2k}\chi^{2k}=\prod^m_{i=1}(\chi^2-a_i^2)$, and
$P(\partial_x,\partial_t)$ be the partial differential operators defined by
$$P(\partial_x,\partial_t):=\sum^m_{k=0}b_{2k}P(\partial_x)^{m-k}\frac{\partial^{2k}}{\partial{t^{2k}}}=\prod^m_{i=1}\left(\frac{\partial^2}{\partial t^2}-a_i^2P(\partial_x)\right),
\quad{x}\in\mathbb{R}^n,\;t\in\mathbb{R}_+^1,\;m\geq1.$$
Here $P(\partial_x)$ be a partial differential operators of any order. Then we have
\begin{equation}\label{j2}
    \left\{\begin{array}{l@{\qquad}l}\displaystyle
    P(\partial_x,\partial_t)u=f(x,t),&x\in\mathbb{R}^n,\;t\in\mathbb{R}_+^1,\\\displaystyle
    \left.\frac{\partial^ru}{\partial{t^r}}\right|_{t=0}=\varphi_r(x),&r=0,1,2,\ldots,2m-1.
    \end{array}\right.
\end{equation}
\begin{eqnarray}\label{j2'}
  u(x,t) &=&
  \int^t_0\int^{t-\tau}_0\frac{(t-\tau-\tau')^{2m-3}}{(2m-3)!}\sum^m_{j=1}\frac{a_j^{2m-2}}{\prod^m_{i=1 \atop i\neq{j}}(a_j^2-a_i^2)}
\frac{\sinh(\tau'a_jP(\partial_x)^{1/2})}{a_jP(\partial_x)^{1/2}}f(x,\tau)\,d\tau'd\tau\nonumber\\
   & & +\,\sum^m_{k=0}b_{2k}P(\partial_x)^{m-k}\sum^{2k-1}_{r=0}\frac{\partial^{2k-1-r}}{\partial t^{2k-1-r}}\int^t_0
   \frac{(t-\tau)^{2m-3}}{(2m-3)!}\nonumber\\
   & & \times\,\sum^m_{j=1}\frac{a_j^{2m-2}}{\prod^m_{i=1 \atop i\neq{j}}(a_j^2-a_i^2)}
   \frac{\sinh(\tau{a_j}P(\partial_x)^{1/2})}{a_jP(\partial_x)^{1/2}}\,\varphi_r(x)\,d\tau.
\end{eqnarray}

Now let us prove the Theorems 3.3 and Theorems 3.4. Since the Theorem BI5 and Theorem BI6, we just need to prove the following Corollary 3.2 and Corollary 3.3:

\textbf{Corollary 3.2.} Let $a_1,a_2,\ldots,a_m$ be the real or complex roots different from each other for any algebraic equation
$b_0+b_1\chi+b_2\chi^2+\cdots+b_m\chi^m=0$, and
$P(\partial_x,\partial_t)$ be the partial differential operators defined by
$$P(\partial_x,\partial_t):=\sum^m_{k=0}b_kP(\partial_x)^{m-k}\frac{\partial^k}{\partial{t^k}},
\quad{x}\in\mathbb{R}^n,\;t\in\mathbb{R}_+^1,\;m\geq1.$$
Here $P(\partial_x)$ be a partial differential operators of any order. Then we have
\begin{equation}\label{j3}
    \left\{\begin{array}{l@{\qquad}l}\displaystyle
    P(\partial_x,\partial_t)u=0,&x\in\mathbb{R}^n,\;t\in\mathbb{R}_+^1,\\\displaystyle
    \left.\frac{\partial^ru}{\partial{t^r}}\right|_{t=0}=\varphi_r(x),&r=0,1,2,\ldots,m-1.
    \end{array}\right.
\end{equation}
\begin{eqnarray}\label{j3'}
  u(x,t)=\sum^m_{k=0}b_kP(\partial_x)^{m-k}\sum^{k-1}_{r=0}\frac{\partial^{k-1-r}}{\partial
t^{k-1-r}}\int^t_0\frac{(t-\tau)^{m-2}}{(m-2)!}
\sum^m_{j=1}\frac{a_j^{m-1}e^{\tau{a_j}P(\partial_x)}}{\prod^m_{i=1
\atop i\neq{j}}(a_j-a_i)}\,\varphi_r(x)\,d\tau.
\end{eqnarray}

\textbf{Proof.} Taking the Laplace transform of both sides of Eq (\ref{j3}) and considering the initial condition gives
$$\sum^m_{k=0}b_kP(\partial_x)^{m-k}\left(s^kU(x,s)-\sum^{k-1}_{r=0}s^{k-1-r}\varphi_r(x)\right)=0,$$
where $U(x,s)=\mathcal{L}u(x,t)$. Considering $\prod^m_{i=1}(s-a_iP(\partial_x))=\sum^m_{k=0}b_ks^kP(\partial_x)^{m-k}$
we have
\begin{equation}\label{j3bh}
  \prod^m_{i=1}(s-a_iP(\partial_x))U(x,s)-\sum^m_{k=0}b_kP(\partial_x)^{m-k}\sum^{k-1}_{r=0}s^{k-1-r}\varphi_r(x)=0.
\end{equation}

We need to introduce the abstract operators $G_m(P(\partial_x),t)$ defined by
$$G_m(P(\partial_x),t):=\mathcal{L}^{-1}\frac{1}{\prod^m_{i=1}(s-a_iP(\partial_x))}.$$
On the other hand, taking the Laplace transform of both sides of Eq (\ref{I1}) and considering its initial condition gives
\begin{equation}\label{I1bh}
  \prod^m_{i=1}(s-a_iP(\partial_x))U(x,s)=F(x,s)\qquad(F(x,s)=\mathcal{L}f(x,t)).
\end{equation}
By solving $U(x,s)$ from (\ref{I1bh}) and using the convolution theorem, we have its inverse transform:
$$u(x,t)=\mathcal{L}^{-1}U(x,s)=\mathcal{L}^{-1}\frac{1}{\prod^m_{i=1}(s-a_iP(\partial_x))}F(x,s)=G_m(P(\partial_x),t)*f(x,t).$$
By comparing (\ref{I1'}) with $u(x,t)=G_m(P(\partial_x),t)*f(x,t)$, we obtain the expression of the abstract operators $G_m(P(\partial_x),t)$:
\begin{equation}\label{j4}
G_m(P(\partial_x),t)=\int^t_0\frac{(t-\tau)^{m-2}}{(m-2)!}\sum^m_{j=1}\frac{a_j^{m-1}}{\prod^m_{i=1 \atop i\neq{j}}(a_j-a_i)}e^{\tau{a_j}P(\partial_x)}d\tau.
\end{equation}
It is easily seen from (\ref{j4}) that
$$\left.\frac{\partial^k}{\partial t^k}G_m(P(\partial_x),t)\right|_{t=0}=0\qquad(k=0,1,2,\ldots,m-2).$$

By solving $U(x,s)$ from (\ref{j3bh}), we have its inverse transform:
\begin{eqnarray}\label{j3''}
u(x,t)&=&\mathcal{L}^{-1}U(x,s)= \sum^m_{k=0}b_kP(\partial_x)^{m-k}\sum^{k-1}_{r=0}\mathcal{L}^{-1}\frac{s^{k-1-r}}{\prod^m_{i=1}(s-a_iP(\partial_x))}\,\varphi_r(x)\nonumber\\
   &=& \sum^m_{k=0}b_kP(\partial_x)^{m-k}\sum^{k-1}_{r=0}\frac{\partial^{k-1-r}}{\partial{t^{k-1-r}}}\,G_m(P(\partial_x),t)\varphi_r(x).
\end{eqnarray}
Applying (\ref{j4}) to (\ref{j3''}), Corollary 3.2 is proved.

\textbf{Corollary 3.3.} Let $a_1,a_2,\ldots,a_m$ be any positive real numbers different from each other such that
$\sum^m_{k=0}b_{2k}\chi^{2k}=\prod^m_{i=1}(\chi^2-a_i^2)$, and
$P(\partial_x,\partial_t)$ be the partial differential operators defined by
$$P(\partial_x,\partial_t):=\sum^m_{k=0}b_{2k}P(\partial_x)^{m-k}\frac{\partial^{2k}}{\partial{t^{2k}}},
\quad{x}\in\mathbb{R}^n,\;t\in\mathbb{R}_+^1,\;m\geq1.$$
Here $P(\partial_x)$ is a partial differential operators of any order. Then we have
\begin{equation}\label{j5}
    \left\{\begin{array}{l@{\qquad}l}\displaystyle
    P(\partial_x,\partial_t)u=0,&x\in\mathbb{R}^n,\;t\in\mathbb{R}_+^1,\\\displaystyle
    \left.\frac{\partial^ru}{\partial{t^r}}\right|_{t=0}=\varphi_r(x),&r=0,1,2,\ldots,2m-1.
    \end{array}\right.
\end{equation}
\begin{eqnarray}\label{j5'}
  u(x,t) &=&
   \sum^m_{k=0}b_{2k}P(\partial_x)^{m-k}\sum^{2k-1}_{r=0}\frac{\partial^{2k-1-r}}{\partial t^{2k-1-r}}\int^t_0
   \frac{(t-\tau)^{2m-3}}{(2m-3)!}\nonumber\\
   & & \times\,\sum^m_{j=1}\frac{a_j^{2m-2}}{\prod^m_{i=1 \atop i\neq{j}}(a_j^2-a_i^2)}
   \frac{\sinh(\tau{a_j}P(\partial_x)^{1/2})}{a_jP(\partial_x)^{1/2}}\,\varphi_r(x)\,d\tau.
\end{eqnarray}

\textbf{Proof.} Taking the Laplace transform of both sides of Eq (\ref{j5}) and considering the initial
condition gives
$$\sum^m_{k=0}b_{2k}P(\partial_x)^{m-k}\left(s^{2k}U(x,s)-\sum^{2k-1}_{r=0}s^{2k-1-r}\varphi_r(x)\right)=0,$$
where $U(x,s)=\mathcal{L}u(x,t)$. Considering
$$\prod^m_{i=1}(s^2-a_i^2P(\partial_x))=\sum^m_{k=0}b_{2k}s^{2k}P(\partial_x)^{m-k},$$
we have
\begin{equation}\label{j5bh}
  \prod^m_{i=1}(s^2-a_i^2P(\partial_x))U(x,s)-\sum^m_{k=0}b_{2k}P(\partial_x)^{m-k}\sum^{2k-1}_{r=0}s^{2k-1-r}\varphi_r(x)=0.
\end{equation}

We need to introduce the abstract operators $G_m(P(\partial_x),t)$ defined by
$$G_m(P(\partial_x),t):=\mathcal{L}^{-1}\frac{1}{\prod^m_{i=1}(s^2-a_i^2P(\partial_x))}.$$
On the other hand, taking the Laplace transform of both sides of Eq (\ref{I2}) and considering its initial condition gives
\begin{equation}\label{I2bh}
  \prod^m_{i=1}(s^2-a_i^2P(\partial_x))U(x,s)=F(x,s)\qquad(F(x,s)=\mathcal{L}f(x,t)).
\end{equation}
By solving $U(x,s)$ from (\ref{I2bh}) and using the convolution theorem, we have its inverse transform:
$$u(x,t)=\mathcal{L}^{-1}U(x,s)=\mathcal{L}^{-1}\frac{1}{\prod^m_{i=1}(s^2-a_i^2P(\partial_x))}F(x,s)=G_m(P(\partial_x),t)*f(x,t).$$
By comparing (\ref{I2'}) with $u(x,t)=G_m(P(\partial_x),t)*f(x,t)$, we obtain the expression of the abstract operators $G_m(P(\partial_x),t)$:
\begin{equation}\label{j6}
G_m(P(\partial_x),t)=\int^t_0\frac{(t-\tau)^{2m-3}}{(2m-3)!}
\sum^m_{j=1}\frac{a_j^{2m-2}}{\prod^m_{i=1 \atop
i\neq{j}}(a_j^2-a_i^2)}\frac{\sinh(\tau{a_j}P(\partial_x)^{1/2})}{a_jP(\partial_x)^{1/2}}\,d\tau.
\end{equation}
It is easily seen from (\ref{j6}) that
$$\left.\frac{\partial^k}{\partial t^k}G_m(P(\partial_x),t)\right|_{t=0}=0\qquad(k=0,1,2,\ldots,2m-2).$$

By solving $U(x,s)$ from (\ref{j5bh}), we have its inverse transform:
\begin{eqnarray}\label{j5''}
   u(x,t)&=&\mathcal{L}^{-1}U(x,s)= \sum^m_{k=0}b_{2k}P(\partial_x)^{m-k}\sum^{2k-1}_{r=0}\mathcal{L}^{-1}\frac{s^{2k-1-r}}{\prod^m_{i=1}(s^2-a_i^2P(\partial_x))}\,\varphi_r(x)\nonumber\\
&=&\sum^m_{k=0}b_{2k}P(\partial_x)^{m-k}\sum^{2k-1}_{r=0}\frac{\partial^{2k-1-r}}{\partial{t^{2k-1-r}}}\,G_m(P(\partial_x),t)\varphi_r(x).
\end{eqnarray}
Applying (\ref{j6}) to (\ref{j5''}), thus Corollary 3.3 is proved.

Another result was given earlier by Guang-Qing Bi \cite[p. 80, Theorem 1]{bi99}:

\textbf{Theorem BI7.} For an arbitrary order partial differential operators $P(\partial_x)$, we have
\begin{equation}\label{I1bi}
    \left\{\begin{array}{l@{\qquad}l}\displaystyle
    \left(\frac{\partial}{\partial{t}}-P(\partial_x)\right)^mu=f(x,t),&x\in\mathbb{R}^n,\;t\in\mathbb{R}_+^1,\;m\geq1,\\\displaystyle
    \left.\frac{\partial^ru}{\partial{t^r}}\right|_{t=0}=\varphi_r(x),&r=0,1,2,\ldots,m-1.
    \end{array}\right.
\end{equation}
\begin{eqnarray}\label{I1'bi}
u(x,t)&=&\int^t_0\frac{(t-\tau)^{m-1}}{(m-1)!}e^{(t-\tau)P(\partial_x)}f(x,\tau)d\tau\nonumber\\
      & & +\,e^{tP(\partial_x)}\sum^{m-1}_{k=0}\sum^k_{r=0}(-1)^{k-r}{k\choose{r}}\frac{t^k}{k!}P(\partial_x)^{k-r}\varphi_r(x).
\end{eqnarray}

It is easily seen that our new method of solving initial value problems for any linear higher-order partial differential equations is universal, which is more convenient than the traditional Fourier transform method.

\subsection{Analytic solutions of Cauchy problem for $n+1$-dimensional multiple inhomogeneous wave equation}
\noindent

\textbf{Theorem 3.5.} (See G.-Q. Bi and Y.-K. Bi \cite[pp. 514-515, Theorem 2]{bi11}) Let
$\Delta:=\sum^n_{k=1}\partial^2_{x_k}$ be the n-dimensional Laplacian. If $n-2=2\nu+1,\;\nu\in\mathbb{N}_0,\;a>0$, then we have
\begin{equation}\label{36}
    \left\{\begin{array}{l@{\qquad}l}\displaystyle
    \left(\frac{\partial^2}{\partial{t^2}}-a^2\Delta\right)^mu=f(x,t),&x\in\mathbb{R}^n,\;t\in\mathbb{R}_+^1,\;m\geq1,\\\displaystyle
    \left.\frac{\partial^ru}{\partial{t^r}}\right|_{t=0}=\varphi_r(x),&r=0,1,2,\ldots,2m-1.
    \end{array}\right.
\end{equation}
\begin{eqnarray}\label{36'}
\lefteqn{u(x,t)=}\nonumber\\
    & & \int^t_0d\tau\int^{t-\tau}_0d\tau'\frac{\left((t-\tau)^2-\tau'^2\right)^{m-2}\tau'^2}{(2m-2)!!\,(2m-4)!!}
    \underbrace{\int^{\tau'}_0\tau'd\tau'\cdots}_\nu\int^{\tau'}_0\frac{(a^2\Delta)^\nu}{S'_n}
  \int_{S'_n}f(\xi',\tau)\,dS'_n\,\tau'd\tau' \nonumber\\
    & & +\,\frac{1}{(2m-2)!!}\sum^{\nu-1}_{r=0}\int^t_0\frac{(t-\tau)^{2m+2r-1}}{(2m+2r-1)!!\,(2r)!!}(a^2\Delta)^rf(x,\tau)\,d\tau \nonumber\\
    & & +\,\sum^{m-1}_{k=0}(-1)^k{m\choose{k}}(a^2\Delta)^{k+\nu}\sum^{2m-1-2k}_{r=0}\frac{\partial^{2m-1-2k-r}}{\partial{t}^{2m-1-2k-r}}
    \int^t_0d\tau\frac{(t^2-\tau^2)^{m-2}\tau^2}{(2m-2)!!\,(2m-4)!!} \nonumber\\
    & & \times\underbrace{\int^\tau_0\tau{d}\tau\cdots}_\nu\int^\tau_0\frac{1}{S_n}\int_{S_n}\varphi_r(\xi)\,dS_n\,\tau{d}\tau
    +\sum^{m-1}_{k=0}(-1)^k{m\choose{k}}\sum^{2m-1-2k}_{r=0}\nonumber\\
    & & \times\sum^{\nu-1}_{i=0}{m-1+i\choose{i}}\frac{t^{2k+2i+r}}{(2k+2i+r)!}(a^2\Delta)^{k+i}\varphi_r(x).
\end{eqnarray}
 Here $n-2=2\nu+1$, $S'_n:=2(2\pi)^{\nu+1}(a\tau')^{n-1}$, $S_n:=2(2\pi)^{\nu+1}(a\tau)^{n-1}$, and $\xi'\in\mathbb{R}^n$ is the integral variable. The integral is on the hypersphere
$(\xi'_1-x_1)^2+(\xi'_2-x_2)^2+\cdots+(\xi'_n-x_n)^2=(a\tau')^2$, and $dS'_n$ is its surface element. $\xi\in\mathbb{R}^n$ is the integral variable on the hypersphere
$(\xi_1-x_1)^2+(\xi_2-x_2)^2+\cdots+(\xi_n-x_n)^2=(a\tau)^2$, and $dS_n$ is its surface element.

\textbf{Proof.} In Theorem 3.2, let $P(\partial_x):=a^2\Delta$. By applying (\ref{37}) to (\ref{32'}) we have
\begin{eqnarray*}
u(x,t)
    &=& \int^t_0\int^{t-\tau}_0\frac{\left((t-\tau)^2-\tau'^2\right)^{m-2}}{(2m-2)!!\,(2m-4)!!}\nonumber\\
    & & \times\left(\tau'\underbrace{\int^{\tau'}_0\tau'd\tau'\cdots}_\nu\int^{\tau'}_0\frac{(a^2\Delta)^\nu}{S'_n}
  \int_{S'_n}f(\xi',\tau)\,dS'_n\,\tau'd\tau'\right)\tau'd\tau'd\tau \nonumber\\
    & & +\int^t_0\int^{t-\tau}_0\frac{\left((t-\tau)^2-\tau'^2\right)^{m-2}}{(2m-2)!!\,(2m-4)!!}
    \sum^{\nu-1}_{r=0}\frac{\tau'^{2r+1}}{(2r+1)!}(a^2\Delta)^rf(x,\tau)\,\tau'd\tau'd\tau \nonumber\\
    & & +\,\sum^{m-1}_{k=0}(-1)^k{m\choose{k}}(a^2\Delta)^k\sum^{2m-1-2k}_{r=0}\frac{\partial^{2m-1-2k-r}}{\partial{t}^{2m-1-2k-r}}
    \int^t_0\frac{(t^2-\tau^2)^{m-2}\tau}{(2m-2)!!\,(2m-4)!!} \nonumber\\
    & & \times\left(\tau\underbrace{\int^\tau_0\tau{d}\tau\cdots}_\nu\int^\tau_0\frac{(a^2\Delta)^\nu}{S_n}\int_{S_n}\varphi_r(\xi)\,dS_n\,\tau{d}\tau\right)d\tau\nonumber\\
    & & +\,\sum^{m-1}_{k=0}(-1)^k{m\choose{k}}(a^2\Delta)^k\sum^{2m-1-2k}_{r=0}\frac{\partial^{2m-1-2k-r}}{\partial{t}^{2m-1-2k-r}}
    \int^t_0\frac{(t^2-\tau^2)^{m-2}\tau}{(2m-2)!!\,(2m-4)!!} \nonumber\\
    & & \times\sum^{\nu-1}_{i=0}\frac{\tau^{2i+1}}{(2i+1)!}(a^2\Delta)^i\varphi_r(x)d\tau.
\end{eqnarray*}
 Here $n-2=2\nu+1$, $S'_n=2(2\pi)^{\nu+1}(a\tau')^{n-1}$, $S_n=2(2\pi)^{\nu+1}(a\tau)^{n-1}$, and $\xi'\in\mathbb{R}^n$ is the integral variable. The integral is on the hypersphere
$(\xi'_1-x_1)^2+(\xi'_2-x_2)^2+\cdots+(\xi'_n-x_n)^2=(a\tau')^2$, and $dS'_n$ is its surface element. $\xi\in\mathbb{R}^n$ is the integral variable on the hypersphere
$(\xi_1-x_1)^2+(\xi_2-x_2)^2+\cdots+(\xi_n-x_n)^2=(a\tau)^2$, and $dS_n$ is its surface element.

By using (\ref{34}), where
$$\frac{\partial^{2m-1-2k-r}}{\partial{t}^{2m-1-2k-r}}\int^t_0\frac{(t^2-\tau^2)^{m-2}\tau}{(2m-2)!!\,(2m-4)!!}\frac{\tau^{2i+1}}{(2i+1)!}d\tau
={m-1+i\choose{i}}\frac{t^{2k+2i+r}}{(2k+2i+r)!}.$$
Similarly,
$$\int^{t-\tau}_0\frac{\left[(t-\tau)^2-\tau'^2\right]^{m-2}}{(2m-2)!!\,(2m-4)!!}\frac{\tau'^{2r+1}}{(2r+1)!}\tau'd\tau'
=\frac{1}{(2m-2)!!}\frac{(t-\tau)^{2m+2r-1}}{(2m+2r-1)!!\,(2r)!!}.$$
Thus Theorem 3.5 is proved.

\textbf{Theorem 3.6.} Let $a_1,a_2,\ldots,a_m$ be any positive real numbers different from each other such that
$\sum^m_{k=0}b_{2k}\chi^{2k}=\prod^m_{i=1}(\chi^2-a_i^2)$, and $P(\partial_x,\partial_t)$ be the partial differential operators defined by
$$P(\partial_x,\partial_t):=\sum^m_{k=0}b_{2k}\Delta^{m-k}\frac{\partial^{2k}}{\partial{t^{2k}}},
\quad{x}\in\mathbb{R}^n,\;t\in\mathbb{R}_+^1,\;m\geq1.$$
 Here $\Delta:=\sum^n_{k=1}\partial^2_{x_k}$ is the n-dimensional Laplacian. If $n-2=2\nu+1,\;\nu=0,1,2,\cdots$, then we have
\begin{equation}\label{j7}
    \left\{\begin{array}{l@{\qquad}l}\displaystyle
    P(\partial_x,\partial_t)u=f(x,t),&x\in\mathbb{R}^n,\;t\in\mathbb{R}_+^1,\\\displaystyle
    \left.\frac{\partial^ru}{\partial{t^r}}\right|_{t=0}=\varphi_r(x),&r=0,1,2,\ldots,2m-1.
    \end{array}\right.
\end{equation}
\begin{eqnarray}\label{j7'}
  u(x,t) &=&
  \int^t_0d\tau\int^{t-\tau}_0d\tau'\frac{(t-\tau-\tau')^{2m-3}}{(2m-3)!}\sum^m_{j=1}\frac{a_j^{2m-2}}{\prod^m_{i=1 \atop i\neq{j}}(a_j^2-a_i^2)}\nonumber\\
  & &
\times\left(\tau'\underbrace{\int^{\tau'}_0\tau'd\tau'\cdots}_\nu\int^{\tau'}_0\frac{(a_j^{2}\Delta)^{\nu}}{S'_{n,j}}\int_{S'_{n,j}}f(\xi',\tau)\,dS'_{n,j}\,\tau'd\tau'\right)\nonumber\\
  & & +\int^t_0\sum^m_{j=1}\frac{a_j^{2m-2}}{\prod^m_{i=1 \atop i\neq{j}}(a_j^2-a_i^2)}\sum^{\nu-1}_{l=0}\frac{a_j^{2l}\,(t-\tau)^{2l+2m-1}}{(2l+2m-1)!}\Delta^lf(x,\tau)d\tau\nonumber\\
  & & +\,\sum^m_{k=0}b_{2k}\Delta^{m-k}\sum^{2k-1}_{r=0}\frac{\partial^{2k-1-r}}{\partial t^{2k-1-r}}\int^t_0d\tau
  \frac{(t-\tau)^{2m-3}}{(2m-3)!}\sum^m_{j=1}\frac{a_j^{2m-2}}{\prod^m_{i=1 \atop i\neq{j}}(a_j^2-a_i^2)}\nonumber\\
  & &
  \times\left(\tau\underbrace{\int^{\tau}_0\tau{d\tau}\cdots}_\nu\int^{\tau}_0\frac{(a_j^{2}\Delta)^{\nu}}{S_{n,j}}\int_{S_{n,j}}\varphi_r(\xi)\,dS_{n,j}\,\tau{d\tau}\right)\nonumber\\
  & & +\,\sum^m_{k=0}b_{2k}\Delta^{m-k}\sum^{2k-1}_{r=0}\sum^m_{j=1}\frac{a_j^{2m-2}}{\prod^m_{i=1 \atop i\neq{j}}(a_j^2-a_i^2)}\sum^{\nu-1}_{l=0}\frac{a_j^{2l}\,t^{2l+2m+r-2k}}{(2l+2m+r-2k)!}\Delta^l\varphi_r(x).
\end{eqnarray}
Here $S'_{n,j}:=2(2\pi)^{\nu+1}(a_j\tau')^{n-1}$, $S_{n,j}:=2(2\pi)^{\nu+1}(a_j\tau)^{n-1}$, and $\xi'\in\mathbb{R}^n$
is the integral variable. The integral is on the hypersphere
$(\xi'_1-x_1)^2+(\xi'_2-x_2)^2+\cdots+(\xi'_n-x_n)^2=(a_j\tau')^2$,
and $dS'_{n,j}$ is its surface element. $\xi\in\mathbb{R}^n$ is the integral variable on the hypersphere
$(\xi_1-x_1)^2+(\xi_2-x_2)^2+\cdots+(\xi_n-x_n)^2=(a_j\tau)^2$, and $dS_{n,j}$ is its surface element.

\textbf{Proof.} In Theorem 3.4, let $P(\partial_x):=\Delta$. By applying (\ref{37}) to (\ref{j2'}) we have
\begin{eqnarray*}
  u(x,t) &=&
  \int^t_0d\tau\int^{t-\tau}_0d\tau'\frac{(t-\tau-\tau')^{2m-3}}{(2m-3)!}\sum^m_{j=1}\frac{a_j^{2m-2}}{\prod^m_{i=1 \atop i\neq{j}}(a_j^2-a_i^2)}\nonumber\\
  & &
  \times\left(\tau'\underbrace{\int^{\tau'}_0\tau'd\tau'\cdots}_\nu\int^{\tau'}_0\frac{(a_j^{2}\Delta)^{\nu}}{S'_{n,j}}\int_{S'_{n,j}}f(\xi',\tau)\,dS'_{n,j}\,\tau'd\tau'
  +\sum^{\nu-1}_{l=0}\frac{a_j^{2l}\tau'^{2l+1}}{(2l+1)!}\Delta^lf(x,\tau)\right)\nonumber\\
   & & +\,\sum^m_{k=0}b_{2k}\Delta^{m-k}\sum^{2k-1}_{r=0}\frac{\partial^{2k-1-r}}{\partial t^{2k-1-r}}\int^t_0d\tau
   \frac{(t-\tau)^{2m-3}}{(2m-3)!}\sum^m_{j=1}\frac{a_j^{2m-2}}{\prod^m_{i=1 \atop i\neq{j}}(a_j^2-a_i^2)}\nonumber\\
   & &
   \times\left(\tau\underbrace{\int^{\tau}_0\tau{d\tau}\cdots}_\nu\int^{\tau}_0\frac{(a_j^{2}\Delta)^{\nu}}{S_{n,j}}\int_{S_{n,j}}\varphi_r(\xi)\,dS_{n,j}\,\tau{d\tau}
   +\sum^{\nu-1}_{l=0}\frac{a_j^{2l}\tau^{2l+1}}{(2l+1)!}\Delta^l\varphi_r(x)\right).
\end{eqnarray*}
Here $S'_{n,j}=2(2\pi)^{\nu+1}(a_j\tau')^{n-1}$, $S_{n,j}=2(2\pi)^{\nu+1}(a_j\tau)^{n-1}$, and $\xi'\in\mathbb{R}^n$
is the integral variable. The integral is on the hypersphere
$(\xi'_1-x_1)^2+(\xi'_2-x_2)^2+\cdots+(\xi'_n-x_n)^2=(a_j\tau')^2$,
and $dS'_{n,j}$ is its surface element. $\xi\in\mathbb{R}^n$ is the integral variable on the hypersphere
$(\xi_1-x_1)^2+(\xi_2-x_2)^2+\cdots+(\xi_n-x_n)^2=(a_j\tau)^2$, and $dS_{n,j}$ is its surface element.
Then Theorem 3.6 is proved.

Anker, Pierfelice and Vallarino \cite{jean} studied the dispersive properties of the wave equation associated with the shifted Laplace-Beltrami operator on real hyperbolic spaces $\mathbb{H}^n$ lately. $\mathbb{H}^n$ can be realized as the symmetric space $G/K$, where $G=\mbox{SO}(1,n)_0$ and $K=\mbox{SO}(n)$. In geodesic polar coordinates on $\mathbb{H}^n$, the Riemannian volume writes (See \cite[pp. 5615-5618]{jean})
\[dx=\mbox{const.}(\sinh r)^{n-1}drd\sigma\]
and the Laplace-Beltrami operator
\[\Delta_{\mathbb{H}^n}=\partial_r^2+(n-1)\coth r\partial_r+\sinh^{-2}r\Delta_{\mathbb{S}^{n-1}}.\]
The spherical functions $\varphi_\lambda$ on $\mathbb{H}^n$ are normalized radial eigenfunctions of $\Delta_{\mathbb{H}^n}$:
\[\left\{\begin{array}{l@{\qquad}l}
\Delta_{\mathbb{H}^n}\varphi_\lambda=-(\lambda^2+\rho^2)\varphi_\lambda,\\
\varphi_\lambda(0)=1,\end{array}\right.\]
where $\lambda\in\mathbb{C}$ and $\rho=(n-1)/2$. They can be expressed in terms of special functions:
\[\varphi_\lambda(r)=\phi_\lambda^{(\frac{n}{2}-1,-\frac{1}{2})}(r)
=\,_2F_1\left(\frac{\rho}{2}+i\frac{\lambda}{2},\frac{\rho}{2}-i\frac{\lambda}{2};\frac{n}{2};-\sinh^2r\right),\]
where $\phi_\lambda^{(\alpha,\beta)}$ denotes the Jacobi functions and $_2F_1$ the Gauss hypergeometric function.

It is easily seen that the partial differential operators $P(\partial_x)$ can also be generalized from constant coefficient to variable coefficient for Theorem BI7 and Theorem 3.2 to 3.4 such as $P(\partial_x):=\Delta_{\mathbb{H}^n}$ or $\Delta_{\mathbb{H}^n}+\rho^2$.
For example, by applying Theorem 3.2 with $m=1$ and $P(\partial_x)=\Delta_{\mathbb{H}^n}+\rho^2$, we obtain the Duhamel's formula, expressed usually by
\[u(t,x)=(\cos tD_x)f(x)+\frac{\sin tD_x}{D_x}g(x)+\int_0^t\frac{\sin(t-s)D_x}{D_x}F(s,x)ds,\]
where $D_x:=\sqrt{-\Delta_{\mathbb{H}^n}-\rho^2}=i\sqrt{\Delta_{\mathbb{H}^n}+\rho^2}$, which gives the solution of the following inhomogeneous linear wave equation on $\mathbb{H}^n$:
\[\left\{\begin{array}{l@{\qquad}l}
\partial_t^2u(t,x)-(\Delta_{\mathbb{H}^n}+\rho^2)u(t,x)=F(t,x),\\
u(0,x)=f(x),\\
\partial_t|_{t=0}u(t,x)=g(x).\end{array}\right.\]

More generally, by applying Theorem 3.2 with $m\geq1$ and $P(\partial_x)=\Delta_{\mathbb{H}^n}+\rho^2$, we obtain the following theorem:

\textbf{Theorem 3.7.} Let $m\geq1$, $\rho=(n-1)/2$, $D_x:=\sqrt{-\Delta_{\mathbb{H}^n}-\rho^2}$, $\Delta_{\mathbb{H}^n}$ be the Laplace-Beltrami operator on $L^2(\mathbb{H}^n)$. Then the solution of the multiple inhomogeneous linear wave equation on $\mathbb{H}^n$ is given as follows:
\begin{equation}\label{32h}
    \left\{\begin{array}{l@{\qquad}l}\displaystyle
    \left(\partial_t^2-(\Delta_{\mathbb{H}^n}+\rho^2)\right)^mu(t,x)=f(t,x),&x\in\mathbb{H}^n,\;t\in\mathbb{R}_+^1,\\\displaystyle
    \partial_t^r|_{t=0}u(t,x)=g_r(x),&r=0,1,2,\ldots,2m-1.
    \end{array}\right.
\end{equation}
\begin{eqnarray}\label{32'h}
  u(t,x) &=&
  \int^t_0\int^{t-\tau}_0\frac{\left((t-\tau)^2-\tau'^2\right)^{m-2}}{(2m-2)!!\,(2m-4)!!}\,
  \frac{\sin\tau'D_x}{D_x}\,f(\tau,x)\,\tau'd\tau'\,d\tau \nonumber\\
   & & +\,\sum^{m-1}_{k=0}(-1)^k{m\choose{k}}(\Delta_{\mathbb{H}^n}+\rho^2)^k\sum^{2m-1-2k}_{r=0}\frac{\partial^{2m-1-2k-r}}{\partial t^{2m-1-2k-r}}\int^t_0
   \frac{(t^2-\tau^2)^{m-2}\tau}{(2m-2)!!\,(2m-4)!!}\nonumber\\
   & & \times\,\frac{\sin\tau D_x}{D_x}\,g_r(x)\,d\tau.
\end{eqnarray}

In fact, the base functions $e^{\xi x}$ on $\mathbb{R}^n$ can be defined by the following  eigenfunctions of $\partial_x$:
\[\left\{\begin{array}{l@{\qquad}l}
\partial_xu(x)=\xi u(x) & \xi\in\mathbb{R}^n,\\
u(0)=1.\end{array}\right.\]
This means that the concept of base functions can be generalized from $e^{\xi x}$ on $\mathbb{R}^n$ to $\varphi_\lambda$ on $\mathbb{H}^n$. Therefore,
by making use of the normalized radial eigenfunctions $\varphi_\lambda(r)$ on $\mathbb{H}^n$, we can extend the definition of abstract operators to obtain the abstract operators taking $\Delta_{\mathbb{H}^n}$ as the operator element, denoted by $f(t,\Delta_{\mathbb{H}^n})$, which is defined as \[f(t,\Delta_{\mathbb{H}^n})\varphi_\lambda:=f(t,-(\lambda^2+\rho^2))\varphi_\lambda,\]
where $f(t,-(\lambda^2+\rho^2))$ is the symbols of abstract operators $f(t,\Delta_{\mathbb{H}^n})$. Thus we have
\[\frac{\sinh(t{a_j}\Delta_{\mathbb{H}^n}^{1/2})}{a_j\Delta_{\mathbb{H}^n}^{1/2}}\varphi_\lambda
=\frac{\sin(a_j\sqrt{\lambda^2+\rho^2}\,t)}{a_j\sqrt{\lambda^2+\rho^2}}\varphi_\lambda.\]
Similarly, $f(t,\Delta_{\mathbb{H}^n}+\rho^2)\varphi_\lambda:=f(t,-\lambda^2)\varphi_\lambda$, where $f(t,-\lambda^2)$ is the symbols of $f(t,\Delta_{\mathbb{H}^n}+\rho^2)$.

Under suitable assumptions, the spherical Fourier transform of a bi-$K$-invariant function $f$ on $G$ is defined by
\[\mathcal{H}f(\lambda):=\int_Gf(g)\varphi_\lambda(g)dg\]
and the following inversion formula and Plancherel formula hold:
\[f(x)=\mbox{const.}\int_0^\infty\varphi_\lambda(x)|\textbf{c}(\lambda)|^{-2}(\mathcal{H}f(\lambda))d\lambda,\]
\[\|f\|_{L^2}^2=\mbox{const.}\int_0^\infty|\textbf{c}(\lambda)|^{-2}|\mathcal{H}f(\lambda)|^2d\lambda.\]
Here the Harish-Chandra \textbf{c}-function is given by
\[\textbf{c}(\lambda)=\frac{\Gamma(2\rho)}{\Gamma(\rho)}\frac{\Gamma(i\lambda)}{\Gamma(i\lambda+\rho)}.\]
Therefore, if $P(\partial_x):=-D_x^2=\Delta_{\mathbb{H}^n}+\rho^2$ for the abstract operators $G_m(P(\partial_x),t)$ introduced in Theorem 3.2 and Theorem 3.4, then for the bi-$K$-invariant function $f$ we have
\[G_m(-D_x^2,t)f(x)=G_m(\Delta_{\mathbb{H}^n}+\rho^2,t)f(x)=\mbox{const.}\int_0^\infty G_m(-\lambda^2,t)\varphi_\lambda(x)|\textbf{c}(\lambda)|^{-2}(\mathcal{H}f(\lambda))d\lambda.\]

Similarly, by applying Theorem 3.4 with $m\geq1$ and $P(\partial_x)=\Delta_{\mathbb{H}^n}+\rho^2$, we have

\textbf{Theorem 3.8.} Let $m\geq1$, $\rho=(n-1)/2$, $D_x:=\sqrt{-\Delta_{\mathbb{H}^n}-\rho^2}$. If $a_1,a_2,\ldots,a_m$ are any positive real numbers different from each other such that
$\sum^m_{k=0}b_{2k}\chi^{2k}=\prod^m_{i=1}(\chi^2-a_i^2)$, and
$P(\partial_x,\partial_t)$ is the partial differential operators defined by
$$P(\partial_x,\partial_t):=\sum^m_{k=0}b_{2k}(\Delta_{\mathbb{H}^n}+\rho^2)^{m-k}\partial_t^{2k}
=\prod^m_{i=1}\left(\partial_t^2-a_i^2(\Delta_{\mathbb{H}^n}+\rho^2)\right),
\quad{x}\in\mathbb{H}^n,\;t\in\mathbb{R}_+^1.$$
Here $\Delta_{\mathbb{H}^n}$ is the Laplace-Beltrami operator on $L^2(\mathbb{H}^n)$. Then the solution of the multiple wave equation on $\mathbb{H}^n$ is given as follows:
\begin{equation}\label{j2h}
    \left\{\begin{array}{l@{\qquad}l}\displaystyle
    P(\partial_x,\partial_t)u(t,x)=f(t,x),&x\in\mathbb{H}^n,\;t\in\mathbb{R}_+^1,\\\displaystyle
    \partial_t^r|_{t=0}u(t,x)=g_r(x),&r=0,1,2,\ldots,2m-1.
    \end{array}\right.
\end{equation}
\begin{eqnarray}\label{j2'h}
  u(t,x) &=&
  \int^t_0\int^{t-\tau}_0\frac{(t-\tau-\tau')^{2m-3}}{(2m-3)!}\sum^m_{j=1}\frac{a_j^{2m-2}}{\prod^m_{i=1 \atop i\neq{j}}(a_j^2-a_i^2)}
\frac{\sin(\tau'a_jD_x)}{a_jD_x}f(\tau,x)\,d\tau'd\tau\nonumber\\
   & & +\,\sum^m_{k=0}b_{2k}(\Delta_{\mathbb{H}^n}+\rho^2)^{m-k}\sum^{2k-1}_{r=0}\frac{\partial^{2k-1-r}}{\partial t^{2k-1-r}}\int^t_0
   \frac{(t-\tau)^{2m-3}}{(2m-3)!}\nonumber\\
   & & \times\,\sum^m_{j=1}\frac{a_j^{2m-2}}{\prod^m_{i=1 \atop i\neq{j}}(a_j^2-a_i^2)}
   \frac{\sin(\tau{a_j}D_x)}{a_jD_x}\,g_r(x)\,d\tau.
\end{eqnarray}

\subsection{Further applications}
\noindent

In this section we discuss the solvability of initial-boundary value problem for the linear higher-order partial differential equations.

Clearly, we can attach proper boundary conditions to the initial value problems (\ref{32}), (\ref{j1}), (\ref{j2}) and  (\ref{I1bi}) introduced by Theorem 3.2, Theorem 3.3, Theorem 3.4 and Theorem BI7. In order to obtain the well-posedness of these initial-boundary value problems, the partial differential operators $P(\partial_x)$ must have the eigenfunctions related to boundary conditions such that the known functions $f(x,\tau),\,\varphi_r(x)$ in (\ref{32'}), (\ref{j1'}), (\ref{j2'}) and (\ref{I1'bi}) can be expanded as the infinite series expressed by the eigenfunctions of $P(\partial_x)$. In order to solve the corresponding initial-boundary value problems, we need to solve the eigenvalue problem of $P(\partial_x)$ under given boundary conditions to determine a set of orthogonal functions. For instance, in the solving formulas (\ref{32'}), (\ref{j1'}), (\ref{j2'}) and (\ref{I1'bi}), if
$f(x,\tau),\varphi_r(x)\in{L^2}(\Omega)$, and $P(\partial_x)$ is the second-order linear self-adjoint elliptic operators, namely
\begin{equation}\label{38'}
  P(\partial_x)u:=\sum^n_{i,j=1}\frac{\partial}{\partial{x}_j}\left(a_{ij}(x)\frac{\partial{u}}{\partial{x}_i}\right)+c(x)u,
\quad{x}\in\Omega\subset\mathbb{R}^n,
\end{equation}
then the boundary conditions can be added for the definite solution problems (\ref{32}), (\ref{j1}), (\ref{j2}) and (\ref{I1bi}):
$\overline{B}u|_{\partial\Omega}=0$, representing $u|_{\partial\Omega}=0$ or
\[\left[\sum^n_{i,j=1}a_{ij}(x)\frac{\partial{u}}{\partial{x}_j}\cos\langle\mathbf{a},x_i\rangle+b(x)u\right]_{\partial\Omega}=0.\]
Here $\mathbf{a}$ is the unit outward normal of $\partial\Omega$. Thus this kind of initial-boundary value problems boils down to solving the eigenvalue problem of the first boundary value problem of second-order linear self-adjoint elliptic operators:
\begin{equation}\label{38}
    \left\{\begin{array}{l@{\qquad}l}\displaystyle
    \sum^n_{i,j=1}\frac{\partial}{\partial{x}_j}\left(a_{ij}(x)\frac{\partial{u}}{\partial{x}_i}\right)+c(x)u=
    -\lambda{u},\quad{x}\in\Omega\subset\mathbb{R}^n,\\\displaystyle
    \overline{B}u|_{\partial\Omega}=0.
    \end{array}\right.
\end{equation}

For the eigenvalue problem (\ref{38}), we first recall here the following known results (See, e.g., Wang \cite[pp. 156-157, Theorem 3.18, Theorem 3.19]{ws}):

$\bullet$ Let $\Omega\subset\mathbb{R}^n$ be a bounded open domain, and $\partial\Omega$ be smooth. Let $a_{ij}=a_{ji}$, and there exists $\theta>0$ such that
$$\sum^n_{i,j=1}a_{ij}(x)\xi_i\xi_j\geq\theta|\xi|^2,\quad{x\in\Omega}.$$
For
$a_{ij}\in{C}^1(\overline{\Omega}),\;c(x)\in{C}(\overline{\Omega}),\;b(x)\in{C}(\partial\Omega)$,
then (\ref{38}) has the following countable eigenvalues:
\[0\leq\lambda_1\leq\lambda_2\leq\cdots\leq\lambda_\nu\leq\cdots,\quad\lim_{\nu\rightarrow\infty}\lambda_\nu=\infty\]
(If $(a_{ij})=I$ is a unit matrix, then $\lambda_1=0$ when $b(x)=c(x)=0$. When $b(x)\geq0,\,c(x)\geq0$ and one of them does not identically equal to zero, $\lambda_1>0$) and the corresponding eigenfunctions $e_1(x),e_2(x),\cdots,e_\nu(x),\cdots,$ satisfy
\begin{equation}\label{39}
  \sum^n_{i,j=1}\frac{\partial}{\partial{x}_j}\left(a_{ij}(x)\frac{\partial{e_i}}{\partial{x}_i}\right)+c(x)e_i=-\lambda_i{e_i},\quad
(e_i,\,e_j)=\delta_{ij}
\end{equation}
and $\{e_j(x)\}^\infty_{j=1}$ are complete in $L^2(\Omega)$. Thus for any $f(x)\in{L^2}(\Omega)$, there exists $c_j$ such that
\[\lim_{\nu\rightarrow\infty}\|f-\sum^\nu_{i=1}c_ie_i\|_{L^2(\Omega)}=0.\]

$\bullet$ Let $\Omega\in\mathbb{R}^n$ be a bounded smooth domain. Then for the eigenvalue problem of the Laplace operators
\begin{equation}\label{ws}
    \left\{\begin{array}{l@{\qquad}l}
    \Delta u=-\lambda{u},&x\in\Omega,\\
    u|_{\partial\Omega}=0,
    \end{array}\right.
\end{equation}
an orthogonal system of the Hilbert space $H_0^1(\Omega)$ is composed of its solutions $\{e_j(x)\}^\infty_{j=1}$.

Therefore, for the solving formulas (\ref{32'}), (\ref{j1'}), (\ref{j2'}) and (\ref{I1'bi}), if $f(x,\tau),\varphi_r(x)\in{L^2}(\Omega)$, we also have
\[\lim_{\nu\rightarrow\infty}\|f(x,\tau)-\sum^\nu_{i=1}c_i(\tau)e_i(x)\|_{L^2(\Omega)}=0,\]
\[\lim_{\nu\rightarrow\infty}\|\varphi_r(x)-\sum^\nu_{i=1}c_ie_i(x)\|_{L^2(\Omega)}=0.\]
Clearly, based on (\ref{38'}) and (\ref{39}), the eigenfunctions $\{e_j(x)\}^\infty_{j=1}$ from the eigenvalue problem (\ref{38}) can be called the base functions of Hilbert space. Therefore, the abstract operators
\[e^{t{a_j}P(\partial_x)},\quad\frac{\sinh(ta_jP(\partial_x)^{1/2})}{a_jP(\partial_x)^{1/2}}\quad\mbox{and}\quad\cosh(ta_jP(\partial_x)^{1/2})
=\frac{\partial}{\partial{t}}\,\frac{\sinh(ta_jP(\partial_x)^{1/2})}{a_jP(\partial_x)^{1/2}}\]
defined on the Hilbert space are also called the abstract operators taking $P(\partial_x)$ as the operator element, denoted by $f(t,P(\partial_x))$ : $L^2(\Omega)\rightarrow{L^2(\Omega)}$, and which acts on the base functions $\{e_j(x)\}^\infty_{j=1}$ such that \[f(t,P(\partial_x))e_i(x):=f(t,-\lambda_i)e_i(x),\]
where $\{f(t,-\lambda_i)\}^\infty_{i=1}$ are the symbols of abstract operators $f(t,P(\partial_x))$ on $L^2(\Omega)$.
Thus we obtain
\begin{equation}\label{40}
  e^{\tau{a_j}P(\partial_x)}e_i(x)=e^{-\tau{a_j}\lambda_i}e_i(x)\quad\mbox{and}\quad\frac{\sinh(\tau{a_j}P(\partial_x)^{1/2})}{a_jP(\partial_x)^{1/2}}e_i(x)=
  \frac{\sin(a_j\sqrt{\lambda_i}\,\tau)}{a_j\sqrt{\lambda_i}}e_i(x),
\end{equation}
where $x\in\Omega\subset\mathbb{R}^n,\;j=1,2,\cdots,m$. Thus we have the following results:
\[\lim_{\nu\rightarrow\infty}\|e^{\tau'{a_j}P(\partial_x)}f(x,\tau)-\sum^\nu_{i=1}(f,e_i)e^{-\tau'{a_j}\lambda_i}e_i(x)\|_{L^2(\Omega)}=0;\]
\[\lim_{\nu\rightarrow\infty}\|e^{\tau{a_j}P(\partial_x)}\varphi_r(x)-\sum^\nu_{i=1}(\varphi_r,e_i)e^{-\tau{a_j}\lambda_i}e_i(x)\|_{L^2(\Omega)}=0;\]
\[\lim_{\nu\rightarrow\infty}\left\|\frac{\sinh(\tau'a_jP(\partial_x)^{1/2})}{a_jP(\partial_x)^{1/2}}f(x,\tau)-
\sum^\nu_{i=1}(f,e_i)\frac{\sin(a_j\sqrt{\lambda_i}\,\tau')}{a_j\sqrt{\lambda_i}}e_i(x)\right\|_{L^2(\Omega)}=0;\]
\[\lim_{\nu\rightarrow\infty}\left\|\frac{\sinh(\tau{a_j}P(\partial_x)^{1/2})}{a_jP(\partial_x)^{1/2}}\varphi_r(x)
-\sum^\nu_{i=1}(\varphi_r,e_i)\frac{\sin(a_j\sqrt{\lambda_i}\,\tau)}{a_j\sqrt{\lambda_i}}e_i(x)\right\|_{L^2(\Omega)}=0.\]

In summary, we can attach respectively the boundary condition $\overline{B}u|_{\partial\Omega}=0$ to the linear higher-order partial differential equations (\ref{32}), (\ref{j1}), (\ref{j2}) and (\ref{I1bi}), such that if the orthogonal complete system $\{e_j(x)\}^\infty_{j=1}$ from the eigenvalue problem (\ref{38}) can be solved, we will obtain the explicit solutions of the corresponding initial-boundary value problems respectively.

Yan'an Second School, Yan'an, Shaanxi, PR China

\emph{E-mail address}: guangqingbi@sohu.com

\end{document}